\input amstex
\documentstyle{amsppt}

\magnification1200
\hsize14cm
\vsize19cm

\def\StemAE{14}
\def\SagaAL{13}
\def\PropAA{12}
\def\OkadAI{11}
\def\MacMAA{10}
\def\MacdAC{9}
\def\KupeAA{8}
\def\IsOWAA{7}
\def\IsWaAB{6}
\def\IsWaAA{5}
\def\GeViAB{4}
\def\FuKrAA{3}
\def\CiucAH{2}
\def\CiucAB{1}

\def\AA{1.1}
\def\AB{1.2}

\def\CA{3.1}
\def\CB{3.2}
\def\CC{3.3}
\def\CD{3.4}

\def\DA{4.1}
\def\DB{4.2}
\def\DC{4.3}

\def\EA{5.3}
\def\EB{5.4}
\def\EC{5.7}
\def\ED{5.8}
\def\EE{5.1}
\def\EF{5.2}
\def\EG{5.5}
\def\EH{5.6}

\def\TA{1}
\def\TB{2}
\def\TC{3}
\def\TD{4}
\def\TE{5}

\def\TG{7}
\def\TH{8}
\def\TI{9}
\def\TJ{10}

\TagsOnRight

\catcode`\@=11
\font@\twelverm=cmr10 scaled\magstep1
\font@\twelveit=cmti10 scaled\magstep1
\font@\twelvebf=cmbx10 scaled\magstep1
\font@\twelvei=cmmi10 scaled\magstep1
\font@\twelvesy=cmsy10 scaled\magstep1
\font@\twelveex=cmex10 scaled\magstep1

\newtoks\twelvepoint@
\def\twelvepoint{\normalbaselineskip15\p@
 \abovedisplayskip15\p@ plus3.6\p@ minus10.8\p@
 \belowdisplayskip\abovedisplayskip
 \abovedisplayshortskip\z@ plus3.6\p@
 \belowdisplayshortskip8.4\p@ plus3.6\p@ minus4.8\p@
 \textonlyfont@\rm\twelverm \textonlyfont@\it\twelveit
 \textonlyfont@\sl\twelvesl \textonlyfont@\bf\twelvebf
 \textonlyfont@\smc\twelvesmc \textonlyfont@\tt\twelvett
%
 \ifsyntax@ \def\big##1{{\hbox{$\left##1\right.$}}}%
  \let\Big\big \let\bigg\big \let\Bigg\big
 \else
  \textfont\z@=\twelverm  \scriptfont\z@=\tenrm  \scriptscriptfont\z@=\sevenrm
  \textfont\@ne=\twelvei  \scriptfont\@ne=\teni  \scriptscriptfont\@ne=\seveni
  \textfont\tw@=\twelvesy \scriptfont\tw@=\tensy \scriptscriptfont\tw@=\sevensy
  \textfont\thr@@=\twelveex \scriptfont\thr@@=\tenex
        \scriptscriptfont\thr@@=\tenex
  \textfont\itfam=\twelveit \scriptfont\itfam=\tenit
        \scriptscriptfont\itfam=\tenit
  \textfont\bffam=\twelvebf \scriptfont\bffam=\tenbf
        \scriptscriptfont\bffam=\sevenbf
  \setbox\strutbox\hbox{\vrule height10.2\p@ depth4.2\p@ width\z@}%
  \setbox\strutbox@\hbox{\lower.6\normallineskiplimit\vbox{%
        \kern-\normallineskiplimit\copy\strutbox}}%
 \setbox\z@\vbox{\hbox{$($}\kern\z@}\bigsize@=1.4\ht\z@
 \fi
 \normalbaselines\rm\ex@.2326ex\jot3.6\ex@\the\twelvepoint@}

\font@\fourteenrm=cmr10 scaled\magstep2
\font@\fourteenit=cmti10 scaled\magstep2
\font@\fourteensl=cmsl10 scaled\magstep2
\font@\fourteensmc=cmcsc10 scaled\magstep2
\font@\fourteentt=cmtt10 scaled\magstep2
\font@\fourteenbf=cmbx10 scaled\magstep2
\font@\fourteeni=cmmi10 scaled\magstep2
\font@\fourteensy=cmsy10 scaled\magstep2
\font@\fourteenex=cmex10 scaled\magstep2
\font@\fourteenmsa=msam10 scaled\magstep2
\font@\fourteeneufm=eufm10 scaled\magstep2
\font@\fourteenmsb=msbm10 scaled\magstep2
\newtoks\fourteenpoint@
\def\fourteenpoint{\normalbaselineskip15\p@
 \abovedisplayskip18\p@ plus4.3\p@ minus12.9\p@
 \belowdisplayskip\abovedisplayskip
 \abovedisplayshortskip\z@ plus4.3\p@
 \belowdisplayshortskip10.1\p@ plus4.3\p@ minus5.8\p@
 \textonlyfont@\rm\fourteenrm \textonlyfont@\it\fourteenit
 \textonlyfont@\sl\fourteensl \textonlyfont@\bf\fourteenbf
 \textonlyfont@\smc\fourteensmc \textonlyfont@\tt\fourteentt
%
 \ifsyntax@ \def\big##1{{\hbox{$\left##1\right.$}}}%
  \let\Big\big \let\bigg\big \let\Bigg\big
 \else
  \textfont\z@=\fourteenrm  \scriptfont\z@=\twelverm  \scriptscriptfont\z@=\tenrm
  \textfont\@ne=\fourteeni  \scriptfont\@ne=\twelvei  \scriptscriptfont\@ne=\teni
  \textfont\tw@=\fourteensy \scriptfont\tw@=\twelvesy \scriptscriptfont\tw@=\tensy
  \textfont\thr@@=\fourteenex \scriptfont\thr@@=\twelveex
        \scriptscriptfont\thr@@=\twelveex
  \textfont\itfam=\fourteenit \scriptfont\itfam=\twelveit
        \scriptscriptfont\itfam=\twelveit
  \textfont\bffam=\fourteenbf \scriptfont\bffam=\twelvebf
        \scriptscriptfont\bffam=\tenbf
  \setbox\strutbox\hbox{\vrule height12.2\p@ depth5\p@ width\z@}%
  \setbox\strutbox@\hbox{\lower.72\normallineskiplimit\vbox{%
        \kern-\normallineskiplimit\copy\strutbox}}%
 \setbox\z@\vbox{\hbox{$($}\kern\z@}\bigsize@=1.7\ht\z@
 \fi
 \normalbaselines\rm\ex@.2326ex\jot4.3\ex@\the\fourteenpoint@}

\font@\seventeenrm=cmr10 scaled\magstep3
\font@\seventeenit=cmti10 scaled\magstep3
\font@\seventeensl=cmsl10 scaled\magstep3
\font@\seventeensmc=cmcsc10 scaled\magstep3
\font@\seventeentt=cmtt10 scaled\magstep3
\font@\seventeenbf=cmbx10 scaled\magstep3
\font@\seventeeni=cmmi10 scaled\magstep3
\font@\seventeensy=cmsy10 scaled\magstep3
\font@\seventeenex=cmex10 scaled\magstep3
\font@\seventeenmsa=msam10 scaled\magstep3
\font@\seventeeneufm=eufm10 scaled\magstep3
\font@\seventeenmsb=msbm10 scaled\magstep3
\newtoks\seventeenpoint@
\def\seventeenpoint{\normalbaselineskip18\p@
 \abovedisplayskip21.6\p@ plus5.2\p@ minus15.4\p@
 \belowdisplayskip\abovedisplayskip
 \abovedisplayshortskip\z@ plus5.2\p@
 \belowdisplayshortskip12.1\p@ plus5.2\p@ minus7\p@
 \textonlyfont@\rm\seventeenrm \textonlyfont@\it\seventeenit
 \textonlyfont@\sl\seventeensl \textonlyfont@\bf\seventeenbf
 \textonlyfont@\smc\seventeensmc \textonlyfont@\tt\seventeentt
%
 \ifsyntax@ \def\big##1{{\hbox{$\left##1\right.$}}}%
  \let\Big\big \let\bigg\big \let\Bigg\big
 \else
  \textfont\z@=\seventeenrm  \scriptfont\z@=\fourteenrm  \scriptscriptfont\z@=\twelverm
  \textfont\@ne=\seventeeni  \scriptfont\@ne=\fourteeni  \scriptscriptfont\@ne=\twelvei
  \textfont\tw@=\seventeensy \scriptfont\tw@=\fourteensy \scriptscriptfont\tw@=\twelvesy
  \textfont\thr@@=\seventeenex \scriptfont\thr@@=\fourteenex
        \scriptscriptfont\thr@@=\fourteenex
  \textfont\itfam=\seventeenit \scriptfont\itfam=\fourteenit
        \scriptscriptfont\itfam=\fourteenit
  \textfont\bffam=\seventeenbf \scriptfont\bffam=\fourteenbf
        \scriptscriptfont\bffam=\twelvebf
  \setbox\strutbox\hbox{\vrule height14.6\p@ depth6\p@ width\z@}%
  \setbox\strutbox@\hbox{\lower.86\normallineskiplimit\vbox{%
        \kern-\normallineskiplimit\copy\strutbox}}%
 \setbox\z@\vbox{\hbox{$($}\kern\z@}\bigsize@=2\ht\z@
 \fi
 \normalbaselines\rm\ex@.2326ex\jot5.2\ex@\the\seventeenpoint@}

\catcode`\@=13
\catcode`\@=11
\font\tenln    = line10
\font\tenlnw   = linew10

\newskip\Einheit \Einheit=0.5cm
\newcount\xcoord \newcount\ycoord
\newdimen\xdim \newdimen\ydim \newdimen\PfadD@cke \newdimen\Pfadd@cke

\newcount\@tempcnta
\newcount\@tempcntb

\newdimen\@tempdima
\newdimen\@tempdimb

\newdimen\@wholewidth
\newdimen\@halfwidth

\newcount\@xarg
\newcount\@yarg
\newcount\@yyarg
\newbox\@linechar
\newbox\@tempboxa
\newdimen\@linelen
\newdimen\@clnwd
\newdimen\@clnht

\newif\if@negarg

\def\@whilenoop#1{}
\def\@whiledim#1\do #2{\ifdim #1\relax#2\@iwhiledim{#1\relax#2}\fi}
\def\@iwhiledim#1{\ifdim #1\let\@nextwhile=\@iwhiledim
        \else\let\@nextwhile=\@whilenoop\fi\@nextwhile{#1}}

\def\@whileswnoop#1\fi{}
\def\@whilesw#1\fi#2{#1#2\@iwhilesw{#1#2}\fi\fi}
\def\@iwhilesw#1\fi{#1\let\@nextwhile=\@iwhilesw
         \else\let\@nextwhile=\@whileswnoop\fi\@nextwhile{#1}\fi}

\def\thinlines{\let\@linefnt\tenln \let\@circlefnt\tencirc
  \@wholewidth\fontdimen8\tenln \@halfwidth .5\@wholewidth}
\def\thicklines{\let\@linefnt\tenlnw \let\@circlefnt\tencircw
  \@wholewidth\fontdimen8\tenlnw \@halfwidth .5\@wholewidth}
\thinlines

\PfadD@cke1pt \Pfadd@cke0.5pt
\def\PfadDicke#1{\PfadD@cke#1 \divide\PfadD@cke by2 \Pfadd@cke\PfadD@cke \multiply\PfadD@cke by2}
\long\def\LOOP#1\REPEAT{\def\BODY{#1}\ITERATE}
\def\ITERATE{\BODY \let\next\ITERATE \else\let\next\relax\fi \next}
\let\REPEAT=\fi
\def\Punkt{\hbox{\raise-2pt\hbox to0pt{\hss$\ssize\bullet$\hss}}}
\def\DuennPunkt(#1,#2){\unskip
  \raise#2 \Einheit\hbox to0pt{\hskip#1 \Einheit
          \raise-2.5pt\hbox to0pt{\hss$\bullet$\hss}\hss}}
\def\NormalPunkt(#1,#2){\unskip
  \raise#2 \Einheit\hbox to0pt{\hskip#1 \Einheit
          \raise-3pt\hbox to0pt{\hss\twelvepoint$\bullet$\hss}\hss}}
\def\DickPunkt(#1,#2){\unskip
  \raise#2 \Einheit\hbox to0pt{\hskip#1 \Einheit
          \raise-4pt\hbox to0pt{\hss\fourteenpoint$\bullet$\hss}\hss}}
\def\Kreis(#1,#2){\unskip
  \raise#2 \Einheit\hbox to0pt{\hskip#1 \Einheit
          \raise-4pt\hbox to0pt{\hss\fourteenpoint$\circ$\hss}\hss}}

\def\Line@(#1,#2)#3{\@xarg #1\relax \@yarg #2\relax
\@linelen=#3\Einheit
\ifnum\@xarg =0 \@vline
  \else \ifnum\@yarg =0 \@hline \else \@sline\fi
\fi}

\def\@sline{\ifnum\@xarg< 0 \@negargtrue \@xarg -\@xarg \@yyarg -\@yarg
  \else \@negargfalse \@yyarg \@yarg \fi
\ifnum \@yyarg >0 \@tempcnta\@yyarg \else \@tempcnta -\@yyarg \fi
\ifnum\@tempcnta>6 \@badlinearg\@tempcnta0 \fi
\ifnum\@xarg>6 \@badlinearg\@xarg 1 \fi
\setbox\@linechar\hbox{\@linefnt\@getlinechar(\@xarg,\@yyarg)}%
\ifnum \@yarg >0 \let\@upordown\raise \@clnht\z@
   \else\let\@upordown\lower \@clnht \ht\@linechar\fi
\@clnwd=\wd\@linechar
\if@negarg \hskip -\wd\@linechar \def\@tempa{\hskip -2\wd\@linechar}\else
     \let\@tempa\relax \fi
\@whiledim \@clnwd <\@linelen \do
  {\@upordown\@clnht\copy\@linechar
   \@tempa
   \advance\@clnht \ht\@linechar
   \advance\@clnwd \wd\@linechar}%
\advance\@clnht -\ht\@linechar
\advance\@clnwd -\wd\@linechar
\@tempdima\@linelen\advance\@tempdima -\@clnwd
\@tempdimb\@tempdima\advance\@tempdimb -\wd\@linechar
\if@negarg \hskip -\@tempdimb \else \hskip \@tempdimb \fi
\multiply\@tempdima \@m
\@tempcnta \@tempdima \@tempdima \wd\@linechar \divide\@tempcnta \@tempdima
\@tempdima \ht\@linechar \multiply\@tempdima \@tempcnta
\divide\@tempdima \@m
\advance\@clnht \@tempdima
\ifdim \@linelen <\wd\@linechar
   \hskip \wd\@linechar
  \else\@upordown\@clnht\copy\@linechar\fi}

\def\@hline{\ifnum \@xarg <0 \hskip -\@linelen \fi
\vrule height\Pfadd@cke width \@linelen depth\Pfadd@cke
\ifnum \@xarg <0 \hskip -\@linelen \fi}

\def\@getlinechar(#1,#2){\@tempcnta#1\relax\multiply\@tempcnta 8
\advance\@tempcnta -9 \ifnum #2>0 \advance\@tempcnta #2\relax\else
\advance\@tempcnta -#2\relax\advance\@tempcnta 64 \fi
\char\@tempcnta}

\def\Vektor(#1,#2)#3(#4,#5){\unskip\leavevmode
  \xcoord#4\relax \ycoord#5\relax
      \raise\ycoord \Einheit\hbox to0pt{\hskip\xcoord \Einheit
         \Vector@(#1,#2){#3}\hss}}

\def\Vector@(#1,#2)#3{\@xarg #1\relax \@yarg #2\relax
\@tempcnta \ifnum\@xarg<0 -\@xarg\else\@xarg\fi
\ifnum\@tempcnta<5\relax
\@linelen=#3\Einheit
\ifnum\@xarg =0 \@vvector
  \else \ifnum\@yarg =0 \@hvector \else \@svector\fi
\fi
\else\@badlinearg\fi}

\def\@hvector{\@hline\hbox to 0pt{\@linefnt
\ifnum \@xarg <0 \@getlarrow(1,0)\hss\else
    \hss\@getrarrow(1,0)\fi}}

\def\@vvector{\ifnum \@yarg <0 \@downvector \else \@upvector \fi}

\def\@svector{\@sline
\@tempcnta\@yarg \ifnum\@tempcnta <0 \@tempcnta=-\@tempcnta\fi
\ifnum\@tempcnta <5
  \hskip -\wd\@linechar
  \@upordown\@clnht \hbox{\@linefnt  \if@negarg
  \@getlarrow(\@xarg,\@yyarg) \else \@getrarrow(\@xarg,\@yyarg) \fi}%
\else\@badlinearg\fi}

\def\@upline{\hbox to \z@{\hskip -.5\Pfadd@cke \vrule width \Pfadd@cke
   height \@linelen depth \z@\hss}}

\def\@downline{\hbox to \z@{\hskip -.5\Pfadd@cke \vrule width \Pfadd@cke
   height \z@ depth \@linelen \hss}}

\def\@upvector{\@upline\setbox\@tempboxa\hbox{\@linefnt\char'66}\raise
     \@linelen \hbox to\z@{\lower \ht\@tempboxa\box\@tempboxa\hss}}

\def\@downvector{\@downline\lower \@linelen
      \hbox to \z@{\@linefnt\char'77\hss}}

\def\@getlarrow(#1,#2){\ifnum #2 =\z@ \@tempcnta='33\else
\@tempcnta=#1\relax\multiply\@tempcnta \sixt@@n \advance\@tempcnta
-9 \@tempcntb=#2\relax\multiply\@tempcntb \tw@
\ifnum \@tempcntb >0 \advance\@tempcnta \@tempcntb\relax
\else\advance\@tempcnta -\@tempcntb\advance\@tempcnta 64
\fi\fi\char\@tempcnta}

\def\@getrarrow(#1,#2){\@tempcntb=#2\relax
\ifnum\@tempcntb < 0 \@tempcntb=-\@tempcntb\relax\fi
\ifcase \@tempcntb\relax \@tempcnta='55 \or
\ifnum #1<3 \@tempcnta=#1\relax\multiply\@tempcnta
24 \advance\@tempcnta -6 \else \ifnum #1=3 \@tempcnta=49
\else\@tempcnta=58 \fi\fi\or
\ifnum #1<3 \@tempcnta=#1\relax\multiply\@tempcnta
24 \advance\@tempcnta -3 \else \@tempcnta=51\fi\or
\@tempcnta=#1\relax\multiply\@tempcnta
\sixt@@n \advance\@tempcnta -\tw@ \else
\@tempcnta=#1\relax\multiply\@tempcnta
\sixt@@n \advance\@tempcnta 7 \fi\ifnum #2<0 \advance\@tempcnta 64 \fi
\char\@tempcnta}

\def\Diagonale(#1,#2)#3{\unskip\leavevmode
  \xcoord#1\relax \ycoord#2\relax
      \raise\ycoord \Einheit\hbox to0pt{\hskip\xcoord \Einheit
         \Line@(1,1){#3}\hss}}
\def\AntiDiagonale(#1,#2)#3{\unskip\leavevmode
  \xcoord#1\relax \ycoord#2\relax 
      \raise\ycoord \Einheit\hbox to0pt{\hskip\xcoord \Einheit
         \Line@(1,-1){#3}\hss}}
\def\Pfad(#1,#2),#3\endPfad{\unskip\leavevmode
  \xcoord#1 \ycoord#2 \thicklines\ZeichnePfad#3\endPfad\thinlines}
\def\ZeichnePfad#1{\ifx#1\endPfad\let\next\relax
  \else\let\next\ZeichnePfad
    \ifnum#1=1
      \raise\ycoord \Einheit\hbox to0pt{\hskip\xcoord \Einheit
         \vrule height\Pfadd@cke width1 \Einheit depth\Pfadd@cke\hss}%
      \advance\xcoord by 1
    \else\ifnum#1=2
      \raise\ycoord \Einheit\hbox to0pt{\hskip\xcoord \Einheit
        \hbox{\hskip-\PfadD@cke\vrule height1 \Einheit width\PfadD@cke depth0pt}\hss}%
      \advance\ycoord by 1
    \else\ifnum#1=3
      \raise\ycoord \Einheit\hbox to0pt{\hskip\xcoord \Einheit
         \Line@(1,1){1}\hss}
      \advance\xcoord by 1
      \advance\ycoord by 1
    \else\ifnum#1=4
      \raise\ycoord \Einheit\hbox to0pt{\hskip\xcoord \Einheit
         \Line@(1,-1){1}\hss}
      \advance\xcoord by 1
      \advance\ycoord by -1
    \fi\fi\fi\fi
  \fi\next}
\def\hSSchritt{\leavevmode\raise-.4pt\hbox to0pt{\hss.\hss}\hskip.2\Einheit
  \raise-.4pt\hbox to0pt{\hss.\hss}\hskip.2\Einheit
  \raise-.4pt\hbox to0pt{\hss.\hss}\hskip.2\Einheit
  \raise-.4pt\hbox to0pt{\hss.\hss}\hskip.2\Einheit
  \raise-.4pt\hbox to0pt{\hss.\hss}\hskip.2\Einheit}
\def\vSSchritt{\vbox{\baselineskip.2\Einheit\lineskiplimit0pt
\hbox{.}\hbox{.}\hbox{.}\hbox{.}\hbox{.}}}
\def\DSSchritt{\leavevmode\raise-.4pt\hbox to0pt{%
  \hbox to0pt{\hss.\hss}\hskip.2\Einheit
  \raise.2\Einheit\hbox to0pt{\hss.\hss}\hskip.2\Einheit
  \raise.4\Einheit\hbox to0pt{\hss.\hss}\hskip.2\Einheit
  \raise.6\Einheit\hbox to0pt{\hss.\hss}\hskip.2\Einheit
  \raise.8\Einheit\hbox to0pt{\hss.\hss}\hss}}
\def\dSSchritt{\leavevmode\raise-.4pt\hbox to0pt{%
  \hbox to0pt{\hss.\hss}\hskip.2\Einheit
  \raise-.2\Einheit\hbox to0pt{\hss.\hss}\hskip.2\Einheit
  \raise-.4\Einheit\hbox to0pt{\hss.\hss}\hskip.2\Einheit
  \raise-.6\Einheit\hbox to0pt{\hss.\hss}\hskip.2\Einheit
  \raise-.8\Einheit\hbox to0pt{\hss.\hss}\hss}}
\def\SPfad(#1,#2),#3\endSPfad{\unskip\leavevmode
  \xcoord#1 \ycoord#2 \ZeichneSPfad#3\endSPfad}
\def\ZeichneSPfad#1{\ifx#1\endSPfad\let\next\relax
  \else\let\next\ZeichneSPfad
    \ifnum#1=1
      \raise\ycoord \Einheit\hbox to0pt{\hskip\xcoord \Einheit
         \hSSchritt\hss}%
      \advance\xcoord by 1
    \else\ifnum#1=2
      \raise\ycoord \Einheit\hbox to0pt{\hskip\xcoord \Einheit
        \hbox{\hskip-2pt \vSSchritt}\hss}%
      \advance\ycoord by 1
    \else\ifnum#1=3
      \raise\ycoord \Einheit\hbox to0pt{\hskip\xcoord \Einheit
         \DSSchritt\hss}
      \advance\xcoord by 1
      \advance\ycoord by 1
    \else\ifnum#1=4
      \raise\ycoord \Einheit\hbox to0pt{\hskip\xcoord \Einheit
         \dSSchritt\hss}
      \advance\xcoord by 1
      \advance\ycoord by -1
    \fi\fi\fi\fi
  \fi\next}
\def\Koordinatenachsen(#1,#2){\unskip
 \hbox to0pt{\hskip-.5pt\vrule height#2 \Einheit width.5pt depth1 \Einheit}%
 \hbox to0pt{\hskip-1 \Einheit \xcoord#1 \advance\xcoord by1
    \vrule height0.25pt width\xcoord \Einheit depth0.25pt\hss}}
\def\Koordinatenachsen(#1,#2)(#3,#4){\unskip
 \hbox to0pt{\hskip-.5pt \ycoord-#4 \advance\ycoord by1
    \vrule height#2 \Einheit width.5pt depth\ycoord \Einheit}%
 \hbox to0pt{\hskip-1 \Einheit \hskip#3\Einheit 
    \xcoord#1 \advance\xcoord by1 \advance\xcoord by-#3 
    \vrule height0.25pt width\xcoord \Einheit depth0.25pt\hss}}
\def\Gitter(#1,#2){\unskip \xcoord0 \ycoord0 \leavevmode
  \LOOP\ifnum\ycoord<#2
    \loop\ifnum\xcoord<#1
      \raise\ycoord \Einheit\hbox to0pt{\hskip\xcoord \Einheit\Punkt\hss}%
      \advance\xcoord by1
    \repeat
    \xcoord0
    \advance\ycoord by1
  \REPEAT}
\def\Gitter(#1,#2)(#3,#4){\unskip \xcoord#3 \ycoord#4 \leavevmode
  \LOOP\ifnum\ycoord<#2
    \loop\ifnum\xcoord<#1
      \raise\ycoord \Einheit\hbox to0pt{\hskip\xcoord \Einheit\Punkt\hss}%
      \advance\xcoord by1
    \repeat
    \xcoord#3
    \advance\ycoord by1
  \REPEAT}
\def\Label#1#2(#3,#4){\unskip \xdim#3 \Einheit \ydim#4 \Einheit
  \def\lo{\advance\xdim by-.5 \Einheit \advance\ydim by.5 \Einheit}%
  \def\llo{\advance\xdim by-.25cm \advance\ydim by.5 \Einheit}%
  \def\loo{\advance\xdim by-.5 \Einheit \advance\ydim by.25cm}%
  \def\o{\advance\ydim by.25cm}%
  \def\ro{\advance\xdim by.5 \Einheit \advance\ydim by.5 \Einheit}%
  \def\rro{\advance\xdim by.25cm \advance\ydim by.5 \Einheit}%
  \def\roo{\advance\xdim by.5 \Einheit \advance\ydim by.25cm}%
  \def\l{\advance\xdim by-.30cm}%
  \def\r{\advance\xdim by.30cm}%
  \def\lu{\advance\xdim by-.5 \Einheit \advance\ydim by-.6 \Einheit}%
  \def\llu{\advance\xdim by-.25cm \advance\ydim by-.6 \Einheit}%
  \def\luu{\advance\xdim by-.5 \Einheit \advance\ydim by-.30cm}%
  \def\u{\advance\ydim by-.30cm}%
  \def\ru{\advance\xdim by.5 \Einheit \advance\ydim by-.6 \Einheit}%
  \def\rru{\advance\xdim by.25cm \advance\ydim by-.6 \Einheit}%
  \def\ruu{\advance\xdim by.5 \Einheit \advance\ydim by-.30cm}%
  #1\raise\ydim\hbox to0pt{\hskip\xdim
     \vbox to0pt{\vss\hbox to0pt{\hss$#2$\hss}\vss}\hss}%
}
\catcode`\@=13
\catcode`\@=11
\font\tenln    = line10
\font\tenlnw   = linew10

\newskip\Einheit \Einheit=0.5cm
\newcount\xcoord \newcount\ycoord
\newdimen\xdim \newdimen\ydim \newdimen\PfadD@cke \newdimen\Pfadd@cke

\newcount\@tempcnta
\newcount\@tempcntb

\newdimen\@tempdima
\newdimen\@tempdimb

\newdimen\@wholewidth
\newdimen\@halfwidth

\newcount\@xarg
\newcount\@yarg
\newcount\@yyarg
\newbox\@linechar
\newbox\@tempboxa
\newdimen\@linelen
\newdimen\@clnwd
\newdimen\@clnht

\newif\if@negarg

\def\@whilenoop#1{}
\def\@whiledim#1\do #2{\ifdim #1\relax#2\@iwhiledim{#1\relax#2}\fi}
\def\@iwhiledim#1{\ifdim #1\let\@nextwhile=\@iwhiledim
        \else\let\@nextwhile=\@whilenoop\fi\@nextwhile{#1}}

\def\@whileswnoop#1\fi{}
\def\@whilesw#1\fi#2{#1#2\@iwhilesw{#1#2}\fi\fi}
\def\@iwhilesw#1\fi{#1\let\@nextwhile=\@iwhilesw
         \else\let\@nextwhile=\@whileswnoop\fi\@nextwhile{#1}\fi}

\def\thinlines{\let\@linefnt\tenln \let\@circlefnt\tencirc
  \@wholewidth\fontdimen8\tenln \@halfwidth .5\@wholewidth}
\def\thicklines{\let\@linefnt\tenlnw \let\@circlefnt\tencircw
  \@wholewidth\fontdimen8\tenlnw \@halfwidth .5\@wholewidth}
\thinlines

\PfadD@cke1pt \Pfadd@cke0.5pt
\def\PfadDicke#1{\PfadD@cke#1 \divide\PfadD@cke by2 \Pfadd@cke\PfadD@cke \multiply\PfadD@cke by2}
\long\def\LOOP#1\REPEAT{\def\BODY{#1}\ITERATE}
\def\ITERATE{\BODY \let\next\ITERATE \else\let\next\relax\fi \next}
\let\REPEAT=\fi
\def\Punkt{\hbox{\raise-2pt\hbox to0pt{\hss$\ssize\bullet$\hss}}}
\def\DuennPunkt(#1,#2){\unskip
  \raise#2 \Einheit\hbox to0pt{\hskip#1 \Einheit
          \raise-2.5pt\hbox to0pt{\hss$\bullet$\hss}\hss}}
\def\NormalPunkt(#1,#2){\unskip
  \raise#2 \Einheit\hbox to0pt{\hskip#1 \Einheit
          \raise-3pt\hbox to0pt{\hss\twelvepoint$\bullet$\hss}\hss}}
\def\DickPunkt(#1,#2){\unskip
  \raise#2 \Einheit\hbox to0pt{\hskip#1 \Einheit
          \raise-4pt\hbox to0pt{\hss\fourteenpoint$\bullet$\hss}\hss}}
\def\Kreis(#1,#2){\unskip
  \raise#2 \Einheit\hbox to0pt{\hskip#1 \Einheit
          \raise-4pt\hbox to0pt{\hss\fourteenpoint$\circ$\hss}\hss}}

\def\Line@(#1,#2)#3{\@xarg #1\relax \@yarg #2\relax
\@linelen=#3\Einheit
\ifnum\@xarg =0 \@vline
  \else \ifnum\@yarg =0 \@hline \else \@sline\fi
\fi}

\def\@sline{\ifnum\@xarg< 0 \@negargtrue \@xarg -\@xarg \@yyarg -\@yarg
  \else \@negargfalse \@yyarg \@yarg \fi
\ifnum \@yyarg >0 \@tempcnta\@yyarg \else \@tempcnta -\@yyarg \fi
\ifnum\@tempcnta>6 \@badlinearg\@tempcnta0 \fi
\ifnum\@xarg>6 \@badlinearg\@xarg 1 \fi
\setbox\@linechar\hbox{\@linefnt\@getlinechar(\@xarg,\@yyarg)}%
\ifnum \@yarg >0 \let\@upordown\raise \@clnht\z@
   \else\let\@upordown\lower \@clnht \ht\@linechar\fi
\@clnwd=\wd\@linechar
\if@negarg \hskip -\wd\@linechar \def\@tempa{\hskip -2\wd\@linechar}\else
     \let\@tempa\relax \fi
\@whiledim \@clnwd <\@linelen \do
  {\@upordown\@clnht\copy\@linechar
   \@tempa
   \advance\@clnht \ht\@linechar
   \advance\@clnwd \wd\@linechar}%
\advance\@clnht -\ht\@linechar
\advance\@clnwd -\wd\@linechar
\@tempdima\@linelen\advance\@tempdima -\@clnwd
\@tempdimb\@tempdima\advance\@tempdimb -\wd\@linechar
\if@negarg \hskip -\@tempdimb \else \hskip \@tempdimb \fi
\multiply\@tempdima \@m
\@tempcnta \@tempdima \@tempdima \wd\@linechar \divide\@tempcnta \@tempdima
\@tempdima \ht\@linechar \multiply\@tempdima \@tempcnta
\divide\@tempdima \@m
\advance\@clnht \@tempdima
\ifdim \@linelen <\wd\@linechar
   \hskip \wd\@linechar
  \else\@upordown\@clnht\copy\@linechar\fi}

\def\@hline{\ifnum \@xarg <0 \hskip -\@linelen \fi
\vrule height\Pfadd@cke width \@linelen depth\Pfadd@cke
\ifnum \@xarg <0 \hskip -\@linelen \fi}

\def\@getlinechar(#1,#2){\@tempcnta#1\relax\multiply\@tempcnta 8
\advance\@tempcnta -9 \ifnum #2>0 \advance\@tempcnta #2\relax\else
\advance\@tempcnta -#2\relax\advance\@tempcnta 64 \fi
\char\@tempcnta}

\def\Vektor(#1,#2)#3(#4,#5){\unskip\leavevmode
  \xcoord#4\relax \ycoord#5\relax
      \raise\ycoord \Einheit\hbox to0pt{\hskip\xcoord \Einheit
         \Vector@(#1,#2){#3}\hss}}

\def\Vector@(#1,#2)#3{\@xarg #1\relax \@yarg #2\relax
\@tempcnta \ifnum\@xarg<0 -\@xarg\else\@xarg\fi
\ifnum\@tempcnta<5\relax
\@linelen=#3\Einheit
\ifnum\@xarg =0 \@vvector
  \else \ifnum\@yarg =0 \@hvector \else \@svector\fi
\fi
\else\@badlinearg\fi}

\def\@hvector{\@hline\hbox to 0pt{\@linefnt
\ifnum \@xarg <0 \@getlarrow(1,0)\hss\else
    \hss\@getrarrow(1,0)\fi}}

\def\@vvector{\ifnum \@yarg <0 \@downvector \else \@upvector \fi}

\def\@svector{\@sline
\@tempcnta\@yarg \ifnum\@tempcnta <0 \@tempcnta=-\@tempcnta\fi
\ifnum\@tempcnta <5
  \hskip -\wd\@linechar
  \@upordown\@clnht \hbox{\@linefnt  \if@negarg
  \@getlarrow(\@xarg,\@yyarg) \else \@getrarrow(\@xarg,\@yyarg) \fi}%
\else\@badlinearg\fi}

\def\@upline{\hbox to \z@{\hskip -.5\Pfadd@cke \vrule width \Pfadd@cke
   height \@linelen depth \z@\hss}}

\def\@downline{\hbox to \z@{\hskip -.5\Pfadd@cke \vrule width \Pfadd@cke
   height \z@ depth \@linelen \hss}}

\def\@upvector{\@upline\setbox\@tempboxa\hbox{\@linefnt\char'66}\raise
     \@linelen \hbox to\z@{\lower \ht\@tempboxa\box\@tempboxa\hss}}

\def\@downvector{\@downline\lower \@linelen
      \hbox to \z@{\@linefnt\char'77\hss}}

\def\@getlarrow(#1,#2){\ifnum #2 =\z@ \@tempcnta='33\else
\@tempcnta=#1\relax\multiply\@tempcnta \sixt@@n \advance\@tempcnta
-9 \@tempcntb=#2\relax\multiply\@tempcntb \tw@
\ifnum \@tempcntb >0 \advance\@tempcnta \@tempcntb\relax
\else\advance\@tempcnta -\@tempcntb\advance\@tempcnta 64
\fi\fi\char\@tempcnta}

\def\@getrarrow(#1,#2){\@tempcntb=#2\relax
\ifnum\@tempcntb < 0 \@tempcntb=-\@tempcntb\relax\fi
\ifcase \@tempcntb\relax \@tempcnta='55 \or
\ifnum #1<3 \@tempcnta=#1\relax\multiply\@tempcnta
24 \advance\@tempcnta -6 \else \ifnum #1=3 \@tempcnta=49
\else\@tempcnta=58 \fi\fi\or
\ifnum #1<3 \@tempcnta=#1\relax\multiply\@tempcnta
24 \advance\@tempcnta -3 \else \@tempcnta=51\fi\or
\@tempcnta=#1\relax\multiply\@tempcnta
\sixt@@n \advance\@tempcnta -\tw@ \else
\@tempcnta=#1\relax\multiply\@tempcnta
\sixt@@n \advance\@tempcnta 7 \fi\ifnum #2<0 \advance\@tempcnta 64 \fi
\char\@tempcnta}

\def\Diagonale(#1,#2)#3{\unskip\leavevmode
  \xcoord#1\relax \ycoord#2\relax
      \raise\ycoord \Einheit\hbox to0pt{\hskip\xcoord \Einheit
         \Line@(1,1){#3}\hss}}
\def\AntiDiagonale(#1,#2)#3{\unskip\leavevmode
  \xcoord#1\relax \ycoord#2\relax 
      \raise\ycoord \Einheit\hbox to0pt{\hskip\xcoord \Einheit
         \Line@(1,-1){#3}\hss}}
\def\Pfad(#1,#2),#3\endPfad{\unskip\leavevmode
  \xcoord#1 \ycoord#2 \thicklines\ZeichnePfad#3\endPfad\thinlines}
\def\ZeichnePfad#1{\ifx#1\endPfad\let\next\relax
  \else\let\next\ZeichnePfad
    \ifnum#1=1
      \raise\ycoord \Einheit\hbox to0pt{\hskip\xcoord \Einheit
         \vrule height\Pfadd@cke width1 \Einheit depth\Pfadd@cke\hss}%
      \advance\xcoord by 1
    \else\ifnum#1=2
      \raise\ycoord \Einheit\hbox to0pt{\hskip\xcoord \Einheit
        \hbox{\hskip-\PfadD@cke\vrule height1 \Einheit width\PfadD@cke depth0pt}\hss}%
      \advance\ycoord by 1
    \else\ifnum#1=3
      \raise\ycoord \Einheit\hbox to0pt{\hskip\xcoord \Einheit
         \Line@(1,1){1}\hss}
      \advance\xcoord by 1
      \advance\ycoord by 1
    \else\ifnum#1=4
      \raise\ycoord \Einheit\hbox to0pt{\hskip\xcoord \Einheit
         \Line@(1,-1){1}\hss}
      \advance\xcoord by 1
      \advance\ycoord by -1
    \fi\fi\fi\fi
  \fi\next}
\def\hSSchritt{\leavevmode\raise-.4pt\hbox to0pt{\hss.\hss}\hskip.2\Einheit
  \raise-.4pt\hbox to0pt{\hss.\hss}\hskip.2\Einheit
  \raise-.4pt\hbox to0pt{\hss.\hss}\hskip.2\Einheit
  \raise-.4pt\hbox to0pt{\hss.\hss}\hskip.2\Einheit
  \raise-.4pt\hbox to0pt{\hss.\hss}\hskip.2\Einheit}
\def\vSSchritt{\vbox{\baselineskip.2\Einheit\lineskiplimit0pt
\hbox{.}\hbox{.}\hbox{.}\hbox{.}\hbox{.}}}
\def\DSSchritt{\leavevmode\raise-.4pt\hbox to0pt{%
  \hbox to0pt{\hss.\hss}\hskip.2\Einheit
  \raise.2\Einheit\hbox to0pt{\hss.\hss}\hskip.2\Einheit
  \raise.4\Einheit\hbox to0pt{\hss.\hss}\hskip.2\Einheit
  \raise.6\Einheit\hbox to0pt{\hss.\hss}\hskip.2\Einheit
  \raise.8\Einheit\hbox to0pt{\hss.\hss}\hss}}
\def\dSSchritt{\leavevmode\raise-.4pt\hbox to0pt{%
  \hbox to0pt{\hss.\hss}\hskip.2\Einheit
  \raise-.2\Einheit\hbox to0pt{\hss.\hss}\hskip.2\Einheit
  \raise-.4\Einheit\hbox to0pt{\hss.\hss}\hskip.2\Einheit
  \raise-.6\Einheit\hbox to0pt{\hss.\hss}\hskip.2\Einheit
  \raise-.8\Einheit\hbox to0pt{\hss.\hss}\hss}}
\def\SPfad(#1,#2),#3\endSPfad{\unskip\leavevmode
  \xcoord#1 \ycoord#2 \ZeichneSPfad#3\endSPfad}
\def\ZeichneSPfad#1{\ifx#1\endSPfad\let\next\relax
  \else\let\next\ZeichneSPfad
    \ifnum#1=1
      \raise\ycoord \Einheit\hbox to0pt{\hskip\xcoord \Einheit
         \hSSchritt\hss}%
      \advance\xcoord by 1
    \else\ifnum#1=2
      \raise\ycoord \Einheit\hbox to0pt{\hskip\xcoord \Einheit
        \hbox{\hskip-2pt \vSSchritt}\hss}%
      \advance\ycoord by 1
    \else\ifnum#1=3
      \raise\ycoord \Einheit\hbox to0pt{\hskip\xcoord \Einheit
         \DSSchritt\hss}
      \advance\xcoord by 1
      \advance\ycoord by 1
    \else\ifnum#1=4
      \raise\ycoord \Einheit\hbox to0pt{\hskip\xcoord \Einheit
         \dSSchritt\hss}
      \advance\xcoord by 1
      \advance\ycoord by -1
    \fi\fi\fi\fi
  \fi\next}
\def\Koordinatenachsen(#1,#2){\unskip
 \hbox to0pt{\hskip-.5pt\vrule height#2 \Einheit width.5pt depth1 \Einheit}%
 \hbox to0pt{\hskip-1 \Einheit \xcoord#1 \advance\xcoord by1
    \vrule height0.25pt width\xcoord \Einheit depth0.25pt\hss}}
\def\Koordinatenachsen(#1,#2)(#3,#4){\unskip
 \hbox to0pt{\hskip-.5pt \ycoord-#4 \advance\ycoord by1
    \vrule height#2 \Einheit width.5pt depth\ycoord \Einheit}%
 \hbox to0pt{\hskip-1 \Einheit \hskip#3\Einheit 
    \xcoord#1 \advance\xcoord by1 \advance\xcoord by-#3 
    \vrule height0.25pt width\xcoord \Einheit depth0.25pt\hss}}
\def\Gitter(#1,#2){\unskip \xcoord0 \ycoord0 \leavevmode
  \LOOP\ifnum\ycoord<#2
    \loop\ifnum\xcoord<#1
      \raise\ycoord \Einheit\hbox to0pt{\hskip\xcoord \Einheit\Punkt\hss}%
      \advance\xcoord by1
    \repeat
    \xcoord0
    \advance\ycoord by1
  \REPEAT}
\def\Gitter(#1,#2)(#3,#4){\unskip \xcoord#3 \ycoord#4 \leavevmode
  \LOOP\ifnum\ycoord<#2
    \loop\ifnum\xcoord<#1
      \raise\ycoord \Einheit\hbox to0pt{\hskip\xcoord \Einheit\Punkt\hss}%
      \advance\xcoord by1
    \repeat
    \xcoord#3
    \advance\ycoord by1
  \REPEAT}
\def\Label#1#2(#3,#4){\unskip \xdim#3 \Einheit \ydim#4 \Einheit
  \def\lo{\advance\xdim by-.5 \Einheit \advance\ydim by.5 \Einheit}%
  \def\llo{\advance\xdim by-.25cm \advance\ydim by.5 \Einheit}%
  \def\loo{\advance\xdim by-.5 \Einheit \advance\ydim by.25cm}%
  \def\o{\advance\ydim by.25cm}%
  \def\ro{\advance\xdim by.5 \Einheit \advance\ydim by.5 \Einheit}%
  \def\rro{\advance\xdim by.25cm \advance\ydim by.5 \Einheit}%
  \def\roo{\advance\xdim by.5 \Einheit \advance\ydim by.25cm}%
  \def\l{\advance\xdim by-.30cm}%
  \def\r{\advance\xdim by.30cm}%
  \def\lu{\advance\xdim by-.5 \Einheit \advance\ydim by-.6 \Einheit}%
  \def\llu{\advance\xdim by-.25cm \advance\ydim by-.6 \Einheit}%
  \def\luu{\advance\xdim by-.5 \Einheit \advance\ydim by-.30cm}%
  \def\u{\advance\ydim by-.30cm}%
  \def\ru{\advance\xdim by.5 \Einheit \advance\ydim by-.6 \Einheit}%
  \def\rru{\advance\xdim by.25cm \advance\ydim by-.6 \Einheit}%
  \def\ruu{\advance\xdim by.5 \Einheit \advance\ydim by-.30cm}%
  #1\raise\ydim\hbox to0pt{\hskip\xdim
     \vbox to0pt{\vss\hbox to0pt{\hss$#2$\hss}\vss}\hss}%
}
\catcode`\@=13






\def\setRevDate $#1 #2 #3${#2}
\def\TeXdrawId{\setRevDate $Date: 1995/12/19 16:40:42 $ TeXdraw V2R0}
\chardef\catamp=\the\catcode`\@
\catcode`\@=11
\long
\def\centertexdraw #1{\hbox to \hsize{\hss
\btexdraw #1\etexdraw
\hss}}
\def\btexdraw {\x@pix=0             \y@pix=0
\x@segoffpix=\x@pix  \y@segoffpix=\y@pix
\t@exdrawdef
\setbox\t@xdbox=\vbox\bgroup\offinterlineskip
\global\d@bs=0
\global\t@extonlytrue
\global\p@osinitfalse
\s@avemove \x@pix \y@pix
\m@pendingfalse
\global\p@osinitfalse
\p@athfalse
\the\everytexdraw}
\def\etexdraw {\ift@extonly \else
\t@drclose
\fi
\egroup
\ifdim \wd\t@xdbox>0pt
\t@xderror {TeXdraw box non-zero size,
possible extraneous text}%
\fi
\vbox {\offinterlineskip
\pixtobp \xminpix \l@lxbp  \pixtobp \yminpix \l@lybp
\pixtobp \xmaxpix \u@rxbp  \pixtobp \ymaxpix \u@rybp
\hbox{\t@xdinclude 
[{\l@lxbp},{\l@lybp}][{\u@rxbp},{\u@rybp}]{\p@sfile}}%
\pixtodim \xminpix \t@xpos  \pixtodim \yminpix \t@ypos
\kern \t@ypos
\hbox {\kern -\t@xpos
\box\t@xdbox
\kern \t@xpos}%
\kern -\t@ypos\relax}}
\def\drawdim #1 {\def\d@dim{#1\relax}}
\def\setunitscale #1 {\edef\u@nitsc{#1}%
\realmult \u@nitsc \s@egsc \d@sc}
\def\relunitscale #1 {\realmult {#1}\u@nitsc \u@nitsc
\realmult \u@nitsc \s@egsc \d@sc}
\def\setsegscale #1 {\edef\s@egsc {#1}%
\realmult \u@nitsc \s@egsc \d@sc}
\def\relsegscale #1 {\realmult {#1}\s@egsc \s@egsc
\realmult \u@nitsc \s@egsc \d@sc}
\def\bsegment {\ifp@ath
\f@lushbs
\f@lushmove
\fi
\begingroup
\x@segoffpix=\x@pix
\y@segoffpix=\y@pix
\setsegscale 1
\global\advance \d@bs by 1\relax}
\def\esegment {\endgroup
\ifnum \d@bs=0
\writetx {es}%
\else
\global\advance \d@bs by -1
\fi}
\def\savecurrpos (#1 #2){\getsympos (#1 #2)\a@rgx\a@rgy
\s@etcsn \a@rgx {\the\x@pix}%
\s@etcsn \a@rgy {\the\y@pix}}
\def\savepos (#1 #2)(#3 #4){\getpos (#1 #2)\a@rgx\a@rgy
\coordtopix \a@rgx \t@pixa
\advance \t@pixa by \x@segoffpix
\coordtopix \a@rgy \t@pixb
\advance \t@pixb by \y@segoffpix
\getsympos (#3 #4)\a@rgx\a@rgy
\s@etcsn \a@rgx {\the\t@pixa}%
\s@etcsn \a@rgy {\the\t@pixb}}
\def\linewd #1 {\coordtopix {#1}\t@pixa
\f@lushbs
\writetx {\the\t@pixa\space sl}}
\def\setgray #1 {\f@lushbs
\writetx {#1 sg}}
\def\lpatt (#1){\listtopix (#1)\p@ixlist
\f@lushbs
\writetx {[\p@ixlist] sd}}
\def\lvec (#1 #2){\getpos (#1 #2)\a@rgx\a@rgy
\s@etpospix \a@rgx \a@rgy
\writeps {\the\x@pix\space \the\y@pix\space lv}}
\def\rlvec (#1 #2){\getpos (#1 #2)\a@rgx\a@rgy
\r@elpospix \a@rgx \a@rgy
\writeps {\the\x@pix\space \the\y@pix\space lv}}
\def\move (#1 #2){\getpos (#1 #2)\a@rgx\a@rgy
\s@etpospix \a@rgx \a@rgy
\s@avemove \x@pix \y@pix}
\def\rmove (#1 #2){\getpos (#1 #2)\a@rgx\a@rgy
\r@elpospix \a@rgx \a@rgy
\s@avemove \x@pix \y@pix}
\def\lcir r:#1 {\coordtopix {#1}\t@pixa
\writeps {\the\t@pixa\space cr}%
\r@elupd \t@pixa \t@pixa
\r@elupd {-\t@pixa}{-\t@pixa}}
\def\fcir f:#1 r:#2 {\coordtopix {#2}\t@pixa
\writeps {\the\t@pixa\space #1 fc}%
\r@elupd \t@pixa \t@pixa
\r@elupd {-\t@pixa}{-\t@pixa}}
\def\lellip rx:#1 ry:#2 {\coordtopix {#1}\t@pixa
\coordtopix {#2}\t@pixb
\writeps {\the\t@pixa\space \the\t@pixb\space el}%
\r@elupd \t@pixa \t@pixb
\r@elupd {-\t@pixa}{-\t@pixb}}
\def\fellip f:#1 rx:#2 ry:#3 {\coordtopix {#2}\t@pixa
\coordtopix {#3}\t@pixb
\writeps {\the\t@pixa\space \the\t@pixb\space #1 fe}%
\r@elupd \t@pixa \t@pixb
\r@elupd {-\t@pixa}{-\t@pixb}}
\def\larc r:#1 sd:#2 ed:#3 {\coordtopix {#1}\t@pixa
\writeps {\the\t@pixa\space #2 #3 ar}}
\def\ifill f:#1 {\writeps {#1 fl}}
\def\lfill f:#1 {\writeps {#1 fp}}
\def\htext #1{\def\testit {#1}%
\ifx \testit\l@paren
\let\next=\h@move
\else
\let\next=\h@text
\fi
\next {#1}}
\def\rtext td:#1 #2{\def\testit {#2}%
\ifx \testit\l@paren
\let\next=\r@move
\else
\let\next=\r@text
\fi
\next td:#1 {#2}}

\def\textref h:#1 v:#2 {\ifx #1R%
\edef\l@stuff {\hss}\edef\r@stuff {}%
\else
\ifx #1C%
\edef\l@stuff {\hss}\edef\r@stuff {\hss}%
\else
\edef\l@stuff {}\edef\r@stuff {\hss}%
\fi
\fi
\ifx #2T%
\edef\t@stuff {}\edef\b@stuff {\vss}%
\else
\ifx #2C%
\edef\t@stuff {\vss}\edef\b@stuff {\vss}%
\else
\edef\t@stuff {\vss}\edef\b@stuff {}%
\fi
\fi}
\def\avec (#1 #2){\getpos (#1 #2)\a@rgx\a@rgy
\s@etpospix \a@rgx \a@rgy
\writeps {\the\x@pix\space \the\y@pix\space (\a@type)
\the\a@lenpix\space \the\a@widpix\space av}}
\def\ravec (#1 #2){\getpos (#1 #2)\a@rgx\a@rgy
\r@elpospix \a@rgx \a@rgy
\writeps {\the\x@pix\space \the\y@pix\space (\a@type)
\the\a@lenpix\space \the\a@widpix\space av}}
\def\arrowheadsize l:#1 w:#2 {\coordtopix{#1}\a@lenpix
\coordtopix{#2}\a@widpix}
\def\arrowheadtype t:#1 {\edef\a@type{#1}}
\def\clvec (#1 #2)(#3 #4)(#5 #6)%
{\getpos (#1 #2)\a@rgx\a@rgy
\coordtopix \a@rgx\t@pixa
\advance \t@pixa by \x@segoffpix
\coordtopix \a@rgy\t@pixb
\advance \t@pixb by \y@segoffpix
\getpos (#3 #4)\a@rgx\a@rgy
\coordtopix \a@rgx\t@pixc
\advance \t@pixc by \x@segoffpix
\coordtopix \a@rgy\t@pixd
\advance \t@pixd by \y@segoffpix
\getpos (#5 #6)\a@rgx\a@rgy
\s@etpospix \a@rgx \a@rgy
\writeps {\the\t@pixa\space \the\t@pixb\space
\the\t@pixc\space \the\t@pixd\space
\the\x@pix\space \the\y@pix\space cv}}
\def\drawbb {\bsegment
\drawdim bp
\linewd 0.24
\setunitscale {\p@sfactor}
\writeps {\the\xminpix\space \the\yminpix\space mv}%
\writeps {\the\xminpix\space \the\ymaxpix\space lv}%
\writeps {\the\xmaxpix\space \the\ymaxpix\space lv}%
\writeps {\the\xmaxpix\space \the\yminpix\space lv}%
\writeps {\the\xminpix\space \the\yminpix\space lv}%
\esegment}
\def\getpos (#1 #2)#3#4{\g@etargxy #1 #2 {} \\#3#4%
\c@heckast #3%
\ifa@st
\g@etsympix #3\t@pixa
\advance \t@pixa by -\x@segoffpix
\pixtocoord \t@pixa #3%
\fi
\c@heckast #4%
\ifa@st
\g@etsympix #4\t@pixa
\advance \t@pixa by -\y@segoffpix
\pixtocoord \t@pixa #4%
\fi}
\def\getsympos (#1 #2)#3#4{\g@etargxy #1 #2 {} \\#3#4%
\c@heckast #3%
\ifa@st \else
\t@xderror {TeXdraw: invalid symbolic coordinate}%
\fi
\c@heckast #4%
\ifa@st \else
\t@xderror {TeXdraw: invalid symbolic coordinate}%
\fi}
\def\listtopix (#1)#2{\def #2{}%
\edef\l@ist {#1 }%
\m@oretrue
\loop
\expandafter\g@etitem \l@ist \\\a@rgx\l@ist
\a@pppix \a@rgx #2%
\ifx \l@ist\empty
\m@orefalse
\fi
\ifm@ore
\repeat}
\def\realmult #1#2#3{\dimen0=#1pt
\dimen2=#2\dimen0
\edef #3{\expandafter\c@lean\the\dimen2}}
\def\intdiv #1#2#3{\t@counta=#1
\t@countb=#2
\ifnum \t@countb<0
\t@counta=-\t@counta
\t@countb=-\t@countb
\fi
\t@countd=1
\ifnum \t@counta<0
\t@counta=-\t@counta
\t@countd=-1
\fi
\t@countc=\t@counta  \divide \t@countc by \t@countb
\t@counte=\t@countc  \multiply \t@counte by \t@countb
\advance \t@counta by -\t@counte
\t@counte=-1
\loop
\advance \t@counte by 1
\ifnum \t@counte<16
\multiply \t@countc by 2
\multiply \t@counta by 2
\ifnum \t@counta<\t@countb \else
\advance \t@countc by 1
\advance \t@counta by -\t@countb
\fi
\repeat
\divide \t@countb by 2
\ifnum \t@counta<\t@countb
\advance \t@countc by 1
\fi
\ifnum \t@countd<0
\t@countc=-\t@countc
\fi
\dimen0=\t@countc sp
\edef #3{\expandafter\c@lean\the\dimen0}}
\def\coordtopix #1#2{\dimen0=#1\d@dim
\dimen2=\d@sc\dimen0
\t@counta=\dimen2
\t@countb=\s@ppix
\divide \t@countb by 2
\ifnum \t@counta<0
\advance \t@counta by -\t@countb
\else
\advance \t@counta by \t@countb
\fi
\divide \t@counta by \s@ppix
#2=\t@counta}
\def\pixtocoord #1#2{\t@counta=#1%
\multiply \t@counta by \s@ppix
\dimen0=\d@sc\d@dim
\t@countb=\dimen0
\intdiv \t@counta \t@countb #2}
\def\pixtodim #1#2{\t@countb=#1%
\multiply \t@countb by \s@ppix
#2=\t@countb sp\relax}
\def\pixtobp #1#2{\dimen0=\p@sfactor pt
\t@counta=\dimen0
\multiply \t@counta by #1%
\ifnum \t@counta < 0
\advance \t@counta by -32768
\else
\advance \t@counta by 32768
\fi
\divide \t@counta by 65536
\edef #2{\the\t@counta}}
\newcount\t@counta    \newcount\t@countb
\newcount\t@countc    \newcount\t@countd
\newcount\t@counte
\newcount\t@pixa      \newcount\t@pixb
\newcount\t@pixc      \newcount\t@pixd
\newdimen\t@xpos      \newdimen\t@ypos
\newcount\xminpix      \newcount\xmaxpix
\newcount\yminpix      \newcount\ymaxpix
\newcount\a@lenpix     \newcount\a@widpix
\newcount\x@pix        \newcount\y@pix
\newcount\x@segoffpix  \newcount\y@segoffpix
\newcount\x@savepix    \newcount\y@savepix
\newcount\s@ppix
\newcount\d@bs
\newcount\t@xdnum
\global\t@xdnum=0
\newbox\t@xdbox
\newwrite\drawfile
\newif\ifm@pending
\newif\ifp@ath
\newif\ifa@st
\newif\ifm@ore
\newif \ift@extonly
\newif\ifp@osinit
\newtoks\everytexdraw
\def\l@paren{(}
\def\a@st{*}
\catcode`\%=12
\def\p@b {
\catcode`\%=14
\catcode`\{=12  \catcode`\}=12  \catcode`\u=1 \catcode`\v=2
\def\l@br u{v  \def\r@br u}v
\catcode `\{=1  \catcode`\}=2   \catcode`\u=11 \catcode`\v=11
{\catcode`\p=12 \catcode`\t=12
\gdef\c@lean #1pt{#1}}
\def\sppix#1/#2 {\dimen0=1#2 \s@ppix=\dimen0
\t@counta=#1%
\divide \t@counta by 2
\advance \s@ppix by \t@counta
\divide \s@ppix by #1%
\t@counta=\s@ppix
\multiply \t@counta by 65536
\advance \t@counta by 32891
\divide \t@counta by 65782
\dimen0=\t@counta sp
\edef\p@sfactor {\expandafter\c@lean\the\dimen0}}
\def\g@etargxy #1 #2 #3 #4\\#5#6{\def #5{#1}%
\ifx #5\empty
\g@etargxy #2 #3 #4 \\#5#6
\else
\def #6{#2}%
\def\next {#3}%
\ifx \next\empty \else
\t@xderror {TeXdraw: invalid coordinate}%
\fi
\fi}
\def\c@heckast #1{\expandafter
\c@heckastll #1\\}
\def\c@heckastll #1#2\\{\def\testit {#1}%
\ifx \testit\a@st
\a@sttrue
\else
\a@stfalse
\fi}
\def\g@etsympix #1#2{\expandafter
\ifx \csname #1\endcsname \relax
\t@xderror {TeXdraw: undefined symbolic coordinate}%
\fi
#2=\csname #1\endcsname}
\def\s@etcsn #1#2{\expandafter
\xdef\csname#1\endcsname {#2}}
\def\g@etitem #1 #2\\#3#4{\edef #4{#2}\edef #3{#1}}
\def\a@pppix #1#2{\edef\next {#1}%
\ifx \next\empty \else
\coordtopix {#1}\t@pixa
\ifx #2\empty
\edef #2{\the\t@pixa}%
\else
\edef #2{#2 \the\t@pixa}%
\fi
\fi}
\def\s@etpospix #1#2{\coordtopix {#1}\x@pix
\advance \x@pix by \x@segoffpix
\coordtopix {#2}\y@pix
\advance \y@pix by \y@segoffpix
\u@pdateminmax \x@pix \y@pix}
\def\r@elpospix #1#2{\coordtopix {#1}\t@pixa
\advance \x@pix by \t@pixa
\coordtopix {#2}\t@pixa
\advance \y@pix by \t@pixa
\u@pdateminmax \x@pix \y@pix}
\def\r@elupd #1#2{\t@counta=\x@pix
\advance\t@counta by #1%
\t@countb=\y@pix
\advance\t@countb by #2%
\u@pdateminmax \t@counta \t@countb}
\def\u@pdateminmax #1#2{\ifnum #1>\xmaxpix
\global\xmaxpix=#1%
\fi
\ifnum #1<\xminpix
\global\xminpix=#1%
\fi
\ifnum #2>\ymaxpix
\global\ymaxpix=#2%
\fi
\ifnum #2<\yminpix
\global\yminpix=#2%
\fi}
\def\s@avemove #1#2{\x@savepix=#1\y@savepix=#2%
\m@pendingtrue
\ifp@osinit \else
\global\p@osinittrue
\global\xminpix=\x@savepix \global\yminpix=\y@savepix
\global\xmaxpix=\x@savepix \global\ymaxpix=\y@savepix
\fi}
\def\f@lushmove {\global\p@osinittrue
\ifm@pending
\writetx {\the\x@savepix\space \the\y@savepix\space mv}%
\m@pendingfalse
\p@athfalse
\fi}
\def\f@lushbs {\loop
\ifnum \d@bs>0
\writetx {bs}%
\global\advance \d@bs by -1
\repeat}
\def\h@move #1#2 #3)#4{\move (#2 #3)%
\h@text {#4}}
\def\h@text #1{\pixtodim \x@pix \t@xpos
\pixtodim \y@pix \t@ypos
\vbox to 0pt{\normalbaselines
\t@stuff
\kern -\t@ypos
\hbox to 0pt{\l@stuff
\kern \t@xpos
\hbox {#1}%
\kern -\t@xpos
\r@stuff}%
\kern \t@ypos
\b@stuff\relax}}
\def\r@move td:#1 #2#3 #4)#5{\move (#3 #4)%
\r@text td:#1 {#5}}
\def\r@text td:#1 #2{\vbox to 0pt{\pixtodim \x@pix \t@xpos
\pixtodim \y@pix \t@ypos
\kern -\t@ypos
\hbox to 0pt{\kern \t@xpos
\rottxt {#1}{\z@sb {#2}}%
\hss}%
\vss}}
\def\z@sb #1{\vbox to 0pt{\normalbaselines
\t@stuff
\hbox to 0pt{\l@stuff \hbox {#1}\r@stuff}%
\b@stuff}}
\ifx \rotatebox\@undefined
\def\rottxt #1#2{\bgroup
#2%
\egroup}
\else
\let\rottxt=\rotatebox
\fi
\ifx \t@xderror\@undefined
\let\t@xderror=\errmessage
\fi
\def\t@exdrawdef {\sppix 300/in
\drawdim in
\edef\u@nitsc {1}%
\setsegscale 1
\arrowheadsize l:0.16 w:0.08
\arrowheadtype t:T
\textref h:L v:B }
\ifx \includegraphics\@undefined
\def\t@xdinclude [#1,#2][#3,#4]#5{%
\begingroup
\message {<#5>}%
\leavevmode
\t@counta=-#1%
\t@countb=-#2%
\setbox0=\hbox{%
\includegraphics{#5}}%
\t@ypos=#4 bp%
\advance \t@ypos by -#2 bp%
\t@xpos=#3 bp%
\advance \t@xpos by -#1 bp%
\dp0=0pt \ht0=\t@ypos  \wd0=\t@xpos
\box0%
\endgroup}
\else
\let\t@xdinclude=\includegraphics
\fi
\def\writeps #1{\f@lushbs
\f@lushmove
\p@athtrue
\writetx {#1}}
\def\writetx #1{\ift@extonly
\global\t@extonlyfalse
\t@xdpsfn \p@sfile
\t@dropen \p@sfile
\fi
\w@rps {#1}}
\def\w@rps #1{\immediate\write\drawfile {#1}}
\def\t@xdpsfn #1{%
\global\advance \t@xdnum by 1
\ifnum \t@xdnum<10
\xdef #1{\jobname.ps\the\t@xdnum}
\else
\xdef #1{\jobname.p\the\t@xdnum}
\fi
}
\def\t@dropen #1{%
\immediate\openout\drawfile=#1%
\w@rps {\p@b PS-Adobe-3.0 EPSF-3.0}%
\w@rps {\p@p BoundingBox: (atend)}%
\w@rps {\p@p Title: TeXdraw drawing: #1}%
\w@rps {\p@p Pages: 1}%
\w@rps {\p@p Creator: \TeXdrawId}%
\w@rps {\p@p CreationDate: \the\year/\the\month/\the\day}%
\w@rps {50 dict begin}%
\w@rps {/mv {stroke moveto} def}%
\w@rps {/lv {lineto} def}%
\w@rps {/st {currentpoint stroke moveto} def}%
\w@rps {/sl {st setlinewidth} def}%
\w@rps {/sd {st 0 setdash} def}%
\w@rps {/sg {st setgray} def}%
\w@rps {/bs {gsave} def /es {stroke grestore} def}%
\w@rps {/fl \l@br gsave setgray fill grestore}%
\w@rps    { currentpoint newpath moveto\r@br\space def}%
\w@rps {/fp {gsave setgray fill grestore st} def}%
\w@rps {/cv {curveto} def}%
\w@rps {/cr \l@br gsave currentpoint newpath 3 -1 roll 0 360 arc}%
\w@rps    { stroke grestore\r@br\space def}%
\w@rps {/fc \l@br gsave setgray currentpoint newpath}%
\w@rps    { 3 -1 roll 0 360 arc fill grestore\r@br\space def}%
\w@rps {/ar {gsave currentpoint newpath 5 2 roll arc stroke grestore} def}%
\w@rps {/el \l@br gsave /svm matrix currentmatrix def}%
\w@rps    { currentpoint translate scale newpath 0 0 1 0 360 arc}%
\w@rps    { svm setmatrix stroke grestore\r@br\space def}%
\w@rps {/fe \l@br gsave setgray currentpoint translate scale newpath}%
\w@rps    { 0 0 1 0 360 arc fill grestore\r@br\space def}%
\w@rps {/av \l@br /hhwid exch 2 div def /hlen exch def}%
\w@rps    { /ah exch def /tipy exch def /tipx exch def}%
\w@rps    { currentpoint /taily exch def /tailx exch def}%
\w@rps    { /dx tipx tailx sub def /dy tipy taily sub def}%
\w@rps    { /alen dx dx mul dy dy mul add sqrt def}%
\w@rps    { /blen alen hlen sub def}%
\w@rps    { gsave tailx taily translate dy dx atan rotate}%
\w@rps    { (V) ah ne {blen 0 gt {blen 0 lineto} if} {alen 0 lineto} ifelse}%
\w@rps    { stroke blen hhwid neg moveto alen 0 lineto blen hhwid lineto}%
\w@rps    { (T) ah eq {closepath} if}%
\w@rps    { (W) ah eq {gsave 1 setgray fill grestore closepath} if}%
\w@rps    { (F) ah eq {fill} {stroke} ifelse}%
\w@rps    { grestore tipx tipy moveto\r@br\space def}%
\w@rps {\p@sfactor\space \p@sfactor\space scale}%
\w@rps {1 setlinecap 1 setlinejoin}%
\w@rps {3 setlinewidth [] 0 setdash}%
\w@rps {0 0 moveto}%
}
\def\t@drclose {%
\bgroup
\w@rps {stroke end showpage}%
\w@rps {\p@p Trailer:}%
\pixtobp \xminpix \l@lxbp  \pixtobp \yminpix \l@lybp
\pixtobp \xmaxpix \u@rxbp  \pixtobp \ymaxpix \u@rybp
\w@rps {\p@p BoundingBox: \l@lxbp\space \l@lybp\space
\u@rxbp\space \u@rybp}%
\w@rps {\p@p EOF}%
\egroup
\immediate\closeout\drawfile
}
\catcode`\@=\catamp

\def\ldreieck{\bsegment
  \rlvec(0.866025403784439 .5) \rlvec(0 -1)
  \rlvec(-0.866025403784439 .5)  
  \savepos(0.866025403784439 -.5)(*ex *ey)
        \esegment
  \move(*ex *ey)
        }
\def\rdreieck{\bsegment
  \rlvec(0.866025403784439 -.5) \rlvec(-0.866025403784439 -.5)  \rlvec(0 1)
  \savepos(0 -1)(*ex *ey)
        \esegment
  \move(*ex *ey)
        }
\def\rhombus{\bsegment
  \rlvec(0.866025403784439 .5) \rlvec(0.866025403784439 -.5) 
  \rlvec(-0.866025403784439 -.5)  \rlvec(0 1)        
  \rmove(0 -1)  \rlvec(-0.866025403784439 .5) 
  \savepos(0.866025403784439 -.5)(*ex *ey)
        \esegment
  \move(*ex *ey)
        }
\def\RhombusA{\bsegment
  \rlvec(0.866025403784439 .5) \rlvec(0.866025403784439 -.5) 
  \rlvec(-0.866025403784439 -.5) \rlvec(-0.866025403784439 .5) 
  \savepos(0.866025403784439 -.5)(*ex *ey)
        \esegment
  \move(*ex *ey)
        }
\def\RhombusB{\bsegment
  \rlvec(0.866025403784439 .5) \rlvec(0 -1)
  \rlvec(-0.866025403784439 -.5) \rlvec(0 1) 
  \savepos(0 -1)(*ex *ey)
        \esegment
  \move(*ex *ey)
        }
\def\RhombusC{\bsegment
  \rlvec(0.866025403784439 -.5) \rlvec(0 -1)
  \rlvec(-0.866025403784439 .5) \rlvec(0 1) 
  \savepos(0.866025403784439 -.5)(*ex *ey)
        \esegment
  \move(*ex *ey)
        }
\def\hdSchritt{\bsegment
  \lpatt(.05 .13)
  \rlvec(0.866025403784439 -.5) 
  \savepos(0.866025403784439 -.5)(*ex *ey)
        \esegment
  \move(*ex *ey)
        }
\def\vdSchritt{\bsegment
  \lpatt(.05 .13)
  \rlvec(0 -1) 
  \savepos(0 -1)(*ex *ey)
        \esegment
  \move(*ex *ey)
        }

\def\la{\lambda}
\def\po#1#2{(#1)_#2}
\def\({\left(}
\def\){\right)}
\def\[{\left[}
\def\]{\right]}
\def\fl#1{\left\lfloor#1\right\rfloor}
\def\cl#1{\left\lceil#1\right\rceil}
\def\flp#1{\lfloor#1\rfloor}
\def\clp#1{\lceil#1\rceil}
\redefine\bar{\overline}
\redefine\hat{\widehat}
\redefine\tilde{\widetilde}
\define\Pf{\operatorname{Pf}}
\define\trans{{}^t\!}

\topmatter 
\title The number of rhombus tilings of a ``punctured" hexagon and
the minor summation formula
\endtitle 
\author S.~Okada\footnote"$^\dagger$"{\hbox{Research supported in part by the
MSRI, Berkeley.}} and C.~Krattenthaler$^\dagger$
\endauthor 
\affil 
Graduate School of Polymathematics,
Nagoya University\\
Furo-cho, Chikusa-ku,
Nagoya 464-01, Japan.\\
email: okada\@math.nagoya-u.ac.jp\\\vskip6pt
Institut f\"ur Mathematik der Universit\"at Wien,\\
Strudlhofgasse 4, A-1090 Wien, Austria.\\
e-mail: KRATT\@Pap.Univie.Ac.At\\
WWW: \tt http://radon.mat.univie.ac.at/People/kratt
\endaffil
\address Graduate School of Polymathematics,
Nagoya University,
Furo-cho, Chikusa-ku,
Nagoya 464-01, Japan
\endaddress
\address Institut f\"ur Mathematik der Universit\"at Wien,
Strudlhofgasse 4, A-1090 Wien, Austria.
\endaddress
\abstract 
We compute the number of all rhombus tilings of a hexagon with sides
$a,b+1,c,a+1,b,c+1$, of which the central triangle is removed,
provided $a,b,c$ have the same parity. The
result is $B(\clp{\frac {a} {2}},\clp{\frac {b} {2}},\clp{\frac {c} {2}})
B(\clp{\frac {a+1} {2}},\flp{\frac {b} {2}},\clp{\frac {c} {2}})
B(\clp{\frac {a} {2}},\clp{\frac {b+1} {2}},\flp{\frac {c} {2}})
B(\flp{\frac {a} {2}},\clp{\frac {b} {2}},\clp{\frac {c+1} {2}})$,
where $B(\alpha,\beta,\gamma)$ is the number of plane partitions inside the
$\alpha\times \beta\times \gamma$ box. The proof uses nonintersecting lattice paths
and a new identity for Schur functions, which is proved by means of
the minor summation formula of Ishikawa and Wakayama. 
A symmetric generalization of
this identity is stated as a conjecture.
\endabstract
\endtopmatter
\document

\leftheadtext{S. Okada and C. Krattenthaler}
\rightheadtext{The number of rhombus tilings of a ``punctured"
hexagon}

\subhead 1. Introduction\endsubhead
In recent years, the enumeration of rhombus tilings of various regions
has attracted a lot of interest and was intensively studied, 
mainly because of the observation (see
\cite{\KupeAA}) that the problem of enumerating all rhombus tilings of a hexagon
with sides $a,b,c,a,b,c$ (see Figure~1; throughout the paper by a
rhombus we always mean a rhombus with side lengths 1 and angles of
$60^\circ$ and $120^\circ$) 
is another way of stating the problem of
counting all plane partitions inside an $a\times b\times c$ box. The
latter problem was solved long ago by MacMahon 
\cite{\MacMAA, Sec.~429, $q\to1$, proof in Sec.~494}. Therefore:

\smallskip
{\it The number of all rhombus tilings of a hexagon
with sides $a,b,c,a,b,c$ equals}
$$B(a,b,c)=\prod _{i=1} ^{a}\prod _{j=1} ^{b}\prod _{k=1} ^{c}\frac {i+j+k-1}
{i+j+k-2}\ .\tag\AA$$
(The form of the expression is due to Macdonald.)

\midinsert
\vskip15pt
\vbox{
\centertexdraw{
  \drawdim truecm  \linewd.02
  \rhombus \rhombus \rhombus \rhombus \ldreieck
  \move (-0.866025403784439 -.5)
  \rhombus \rhombus \rhombus \rhombus \rhombus \ldreieck
  \move (-1.732050807568877 -1)
  \rhombus \rhombus \rhombus \rhombus \rhombus \rhombus \ldreieck
  \move (-1.732050807568877 -1)
  \rdreieck
  \rhombus \rhombus \rhombus \rhombus \rhombus \rhombus \ldreieck
  \move (-1.732050807568877 -2)
  \rdreieck
  \rhombus \rhombus \rhombus \rhombus \rhombus \rhombus \ldreieck
  \move (-1.732050807568877 -3)
  \rdreieck
  \rhombus \rhombus \rhombus \rhombus \rhombus \rhombus 
  \move (-1.732050807568877 -4)
  \rdreieck
  \rhombus \rhombus \rhombus \rhombus \rhombus 
  \move (-1.732050807568877 -5)
  \rdreieck
  \rhombus \rhombus \rhombus \rhombus 
\move(8 0)
\bsegment
  \drawdim truecm  \linewd.02
  \rhombus \rhombus \rhombus \rhombus \ldreieck
  \move (-0.866025403784439 -.5)
  \rhombus \rhombus \rhombus \rhombus \rhombus \ldreieck
  \move (-1.732050807568877 -1)
  \rhombus \rhombus \rhombus \rhombus \rhombus \rhombus \ldreieck
  \move (-1.732050807568877 -1)
  \rdreieck
  \rhombus \rhombus \rhombus \rhombus \rhombus \rhombus \ldreieck
  \move (-1.732050807568877 -2)
  \rdreieck
  \rhombus \rhombus \rhombus \rhombus \rhombus \rhombus \ldreieck
  \move (-1.732050807568877 -3)
  \rdreieck
  \rhombus \rhombus \rhombus \rhombus \rhombus \rhombus 
  \move (-1.732050807568877 -4)
  \rdreieck
  \rhombus \rhombus \rhombus \rhombus \rhombus 
  \move (-1.732050807568877 -5)
  \rdreieck
  \rhombus \rhombus \rhombus \rhombus 
  \linewd.08
  \move(0 0)
  \RhombusA \RhombusB \RhombusB 
  \RhombusA \RhombusA \RhombusB \RhombusA \RhombusB \RhombusB
  \move (-0.866025403784439 -.5)
  \RhombusA \RhombusB \RhombusB \RhombusB \RhombusB
  \RhombusA \RhombusA \RhombusB \RhombusA 
  \move (-1.732050807568877 -1)
  \RhombusB \RhombusB \RhombusA \RhombusB \RhombusB \RhombusA
  \RhombusB \RhombusA \RhombusA 
  \move (1.732050807568877 0)
  \RhombusC \RhombusC \RhombusC 
  \move (1.732050807568877 -1)
  \RhombusC \RhombusC \RhombusC 
  \move (3.464101615137755 -3)
  \RhombusC 
  \move (-0.866025403784439 -.5)
  \RhombusC
  \move (-0.866025403784439 -1.5)
  \RhombusC
  \move (0.866025403784439 -2.5)
  \RhombusC \RhombusC 
  \move (0.866025403784439 -3.5)
  \RhombusC \RhombusC \RhombusC 
  \move (2.598076211353316 -5.5)
  \RhombusC 
  \move (0.866025403784439 -5.5)
  \RhombusC 
  \move (-1.732050807568877 -3)
  \RhombusC 
  \move (-1.732050807568877 -4)
  \RhombusC 
  \move (-1.732050807568877 -5)
  \RhombusC \RhombusC 
\esegment
\htext (-1.5 -9){\eightpoint a. A hexagon with sides $a,b,c,a,b,c$,}
\htext (-1.5 -9.5){\eightpoint \hphantom{a. }where $a=3$, $b=4$, $c=5$}
\htext (6.8 -9){\eightpoint b. A rhombus tiling of a hexagon}
\htext (6.8 -9.5){\eightpoint \hphantom{b. }with sides $a,b,c,a,b,c$}
\rtext td:0 (4.3 -4){$\sideset {} \and c\to 
    {\left.\vbox{\vskip2.2cm}\right\}}$}
\rtext td:60 (2.6 -.5){$\sideset {} \and {} \to 
    {\left.\vbox{\vskip1.7cm}\right\}}$}
\rtext td:120 (-.44 -.25){$\sideset {} \and {} \to 
    {\left.\vbox{\vskip1.3cm}\right\}}$}
\rtext td:0 (-2.3 -3.6){$\sideset {c} \and {}\to 
    {\left\{\vbox{\vskip2.2cm}\right.}$}
\rtext td:240 (0 -7){$\sideset {} \and {} \to 
    {\left.\vbox{\vskip1.7cm}\right\}}$}
\rtext td:300 (3.03 -7.25){$\sideset {} \and {} \to 
    {\left.\vbox{\vskip1.4cm}\right\}}$}
\htext (-.9 0.1){$a$}
\htext (2.8 -.1){$b$}
\htext (3.2 -7.8){$a$}
\htext (-0.3 -7.65){$b$}
}
\centerline{\eightpoint Figure 1}
}
\vskip10pt
\endinsert

In his preprint \cite{\PropAA}, Propp proposes several variations of
this enumeration problem, one of which (Problem~2) asks for the
enumeration of all rhombus tilings of a hexagon with sides
$n,n+1,n,n+1,n,n+1$, of which the central triangle is removed (a
``punctured hexagon"). At this
point, this may seem to be a somewhat artificial problem. But sure
enough, soon after, an even more general problem, namely the problem
of enumerating all rhombus tilings of a hexagon with sides
$a,b+1,c,a+1,b,c+1$, of which the central triangle is removed (see
Figure~2), occured
in work of Kuperberg related to the Penrose impossible triangle [private
communication].

\midinsert

\vskip15pt
\vbox{
\centertexdraw{
  \drawdim truecm  \linewd.02
  \rhombus \rhombus 
  \rhombus \rhombus \rhombus \rhombus \ldreieck
  \move (-0.866025403784439 -.5)
  \rhombus \rhombus \rhombus 
  \rhombus \rhombus \rhombus \rhombus \ldreieck
  \move (-1.732050807568877 -1)
  \rhombus \rhombus \rhombus \rhombus 
  \rhombus \rhombus \rhombus \rhombus \ldreieck
  \move (-1.732050807568877 -1)
  \rdreieck
  \rhombus \rhombus \rhombus \rhombus 
  \rhombus \rhombus \rhombus \rhombus \ldreieck
  \move (-1.732050807568877 -2)
  \rdreieck
  \rhombus \rhombus \rhombus \rhombus 
  \rhombus \rhombus \rhombus \rhombus \ldreieck
  \move (-1.732050807568877 -3)
  \rdreieck
  \rhombus \rhombus \rhombus \rhombus 
  \rhombus \rhombus \rhombus \rhombus 
  \move (-1.732050807568877 -4)
  \rdreieck
  \rhombus \rhombus \rhombus \rhombus 
  \rhombus \rhombus \rhombus 
  \move (-1.732050807568877 -5)
  \rdreieck
  \rhombus \rhombus \rhombus \rhombus 
  \rhombus \rhombus 
  \move (-1.732050807568877 -6)
  \rdreieck
  \rhombus \rhombus \rhombus \rhombus 
  \rhombus 
  \linewd.08
  \move(0 0)
  \RhombusA \RhombusB \RhombusB 
  \RhombusA \RhombusA \RhombusB \RhombusA \RhombusA \RhombusA
  \RhombusB \RhombusB
  \move (-0.866025403784439 -.5)
  \RhombusA \RhombusB \RhombusB \RhombusB \RhombusB
  \RhombusA \RhombusA \RhombusB \RhombusA \RhombusA \RhombusB
  \move (-1.732050807568877 -1)
  \RhombusB \RhombusB \RhombusA \RhombusB \RhombusB \RhombusB
  \RhombusA \RhombusA \RhombusA \RhombusB \RhombusA
  \move (1.732050807568877 -4)
    \bsegment
  \rlvec(0.866025403784439 -.5) \rlvec(-0.866025403784439 -.5) 
  \lfill f:0
  \savepos(0 -1)(*ex *ey)
    \esegment
  \move(*ex *ey)
  \RhombusA \RhombusA \RhombusB \RhombusA \RhombusB
  \move (1.732050807568877 0)
  \RhombusC \RhombusC \RhombusC \RhombusC \RhombusC 
  \move (1.732050807568877 -1)
  \RhombusC \RhombusC \RhombusC \RhombusC \RhombusC 
  \move (3.464101615137755 -3)
  \RhombusC \RhombusC \RhombusC 
  \move (-0.866025403784439 -.5)
  \RhombusC
  \move (-0.866025403784439 -1.5)
  \RhombusC
  \move (0.866025403784439 -2.5)
  \RhombusC \RhombusC 
  \move (0.866025403784439 -3.5)
  \RhombusC 
  \move (0 -5)
  \RhombusC \RhombusC 
  \move (2.598076211353316 -7.5)
  \RhombusC 
  \move (-1.732050807568877 -3)
  \RhombusC 
  \move (-1.732050807568877 -4)
  \RhombusC 
  \move (-1.732050807568877 -5)
  \RhombusC 
  \move (-1.732050807568877 -6)
  \RhombusC \RhombusC \RhombusC \RhombusC
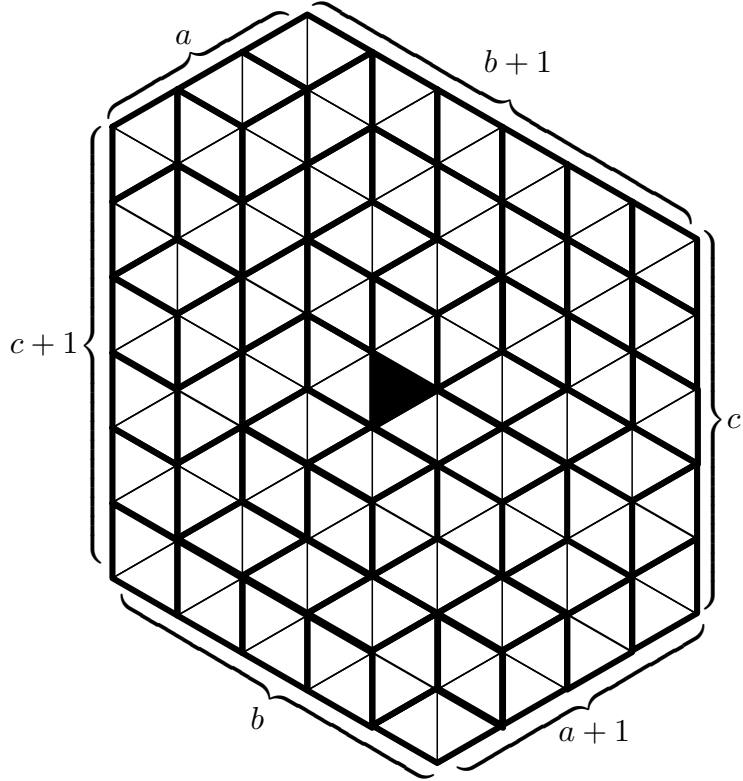
\htext (-1.2 -10.5){\eightpoint A rhombus tiling of a ``punctured" hexagon}
\htext (-2.5 -11){\eightpoint with sides $a,b+1,c,a+1,b,c+1$,
            where $a=3$, $b=5$, $c=5$}
\rtext td:0 (6 -5){$\sideset {} \and c\to 
    {\left.\vbox{\vskip2.2cm}\right\}}$}
\rtext td:60 (3.46 -1){$\sideset {} \and {} \to 
    {\left.\vbox{\vskip2.6cm}\right\}}$}
\rtext td:120 (-.44 -.25){$\sideset {} \and {} \to 
    {\left.\vbox{\vskip1.3cm}\right\}}$}
\rtext td:0 (-3.1 -4){$\sideset {c+1} \and {}\to 
    {\left\{\vbox{\vskip2.5cm}\right.}$}
\rtext td:240 (0.44 -8.25){$\sideset {} \and {} \to 
    {\left.\vbox{\vskip2.1cm}\right\}}$}
\rtext td:300 (4.33 -8.5){$\sideset {} \and {} \to 
    {\left.\vbox{\vskip1.7cm}\right\}}$}
\htext (-.9 0.1){$a$}
\htext (3.2 -.3){$b+1$}
\htext (4.2 -9.2){$a+1$}
\htext (0.1 -9){$b$}
}
\vskip4pt
\centerline{\eightpoint Figure 2}
}
\vskip10pt
\endinsert

The purpose of this paper is to solve this enumeration problem. (We
want to mention that Ciucu \cite{\CiucAH} has an independent
solution of the $b=c$ case of the problem, which builds upon his
matchings factorization theorem \cite{\CiucAB}.) Our
result is the following:
\proclaim{Theorem~\TA}Let $a,b,c$ be positive integers, all of the
same parity. Then the number of all rhombus tilings of a hexagon with sides
$a,b+1,c,a+1,b,c+1$, of which the central triangle is removed, equals
$$\multline
B\(\cl{\frac {\vphantom{b}a} {2}},\cl{\frac {b} {2}},\cl{\frac {\vphantom{b}c} {2}}\)
B\(\cl{\frac {a+1} {2}},\fl{\frac {b} {2}},\cl{\frac {\vphantom{b}c} {2}}\)\\
\times
B\(\cl{\frac {\vphantom{b}a} {2}},\cl{\frac {b+1} {2}},\fl{\frac 
{\vphantom{b}c} {2}}\)
B\(\fl{\frac {\vphantom{b}a} {2}},\cl{\frac {b} {2}},\cl{\frac {c+1}
{2}}\),
\endmultline\tag\AB$$
where $B(\alpha,\beta,\gamma)$ is the number of all plane partitions inside the
$\alpha\times \beta\times \gamma$ box, which is given by \rm(\AA).
\endproclaim

We prove this Theorem by first converting the tiling problem into an
enumeration problem for nonintersecting lattice paths (see Section~2),
deriving a certain summation for the desired number (see Section~3,
Proposition~\TB), and
then evaluating the sum by proving (in Section~4) 
an actually much more general
identity for Schur functions (Theorem~\TC). The summation that we are
interested in then follows immediately by setting all variables
equal to 1. In order to prove this
Schur function identity, which we do in Section~5, 
we make essential use of the minor summation
formula of Ishikawa and Wakayama \cite{\IsWaAA} (see Theorem~\TG).
In fact, it appears that a symmetric generalization of this Schur
function identity (Conjecture~\TE\ in Section~4) holds. 
However, so far we were not
able to establish this identity. Also in Section~4, 
we add a further enumeration
result, Theorem~\TD, on rhombus tilings of a ``punctured" hexagon, which
follows from a different specialization of Theorem~\TC.

\smallskip
In conclusion of the Introduction,
we wish to point out a few interesting features of our result in
Theorem~\TA. 

(1) First of
all, it is another application of the powerful minor summation
formula of Ishikawa and Wakayama. For other applications see
\cite{\IsWaAB, \IsOWAA, \OkadAI}. 

(2) A very striking fact is that
this problem, in the formulation of nonintersecting lattice paths, is
the first instance of existence of a closed form enumeration, despite
the fact that the starting points of the paths are not in the ``right"
order. That is, the ``compatability" condition for the starting and
end points in the main theorem on
nonintersecting lattice paths \cite{\GeViAB, Cor.~2; \StemAE,
Theorem~1.2} is violated (compare Figure~3.c) and the main theorem 
does therefore not apply. (The ``compatability" condition is the
requirement that whenever we consider starting points $A_i$, $A_j$,
with $i<j$, and end points $E_k$, $E_l$, with $k<l$, then any path
from $A_i$ to $E_l$ has to touch any path from $A_j$ to $E_k$.) This is the
reason that we had to resort to something else, which turned out to
be the minor summation formula. 

(3) The result (\AB) is in a
very appealing combinatorial form. It is natural to ask for a
bijective proof of the formula, i.e., for setting up a one-to-one
correspondence between the rhombus tilings in question and a
quadruple of plane partitions as described by the product in (\AB).
However, this seems to be a very difficult problem. In particular,
how does one split a rhombus tiling of a {\it hexagon} into {\it four\/}
objects? 

(4) In the course of our investigations, 
we discovered the symmetric generalization, Conjecture~\TE, of
Theorem~\TC\ which we
mentioned before. Although the minor summation formula still applies,
we were not able to carry out the subsequent steps, i.e., to find the
appropriate generalizations of Lemmas~\TI\ and \TJ.

\subhead 2. From rhombus tilings to nonintersecting lattice paths\endsubhead
In this section we translate our problem of enumerating rhombus
tilings of a ``punctured" hexagon into the language of
nonintersecting lattice paths. Throughout this section we assume that
$a,b,c$ are positive integers, all of the same parity.

We use a slight (but obvious) modification of the standard
translation from rhombus tilings of a hexagon to
nonintersecting lattice paths.
Given a rhombus tiling of a hexagon with sides
$a,b+1,c,a+1,b,c+1$, of which the central triangle is removed (see
Figure~2), we mark the centres of the $a$ edges along the side of
length $a$, we mark the centres of the $a+1$ edges along the side of
length $a+1$, and we mark the centre of the edge of the removed
triangle that is parallel to the sides of respective lengths 
$a$ and $a+1$ (see
Figure~3.a; the marked points are indicated by circles). Starting
from each of the marked points on the side of length $a$, we form a
lattice path by connecting the marked point with the centre of the
edge opposite to it, the latter centre with the centre of the edge
opposite to it, etc\. (see Figure~3.a,b; the paths are indicated by
broken lines). For convenience, we deform the picture so that the
slanted edges of the paths become horizontal (see Figure~3.c). Thus
we obtain a family of lattice paths in the integer lattice $\Bbb
Z^2$, consisting of horizontal unit steps in the positive direction 
and vertical unit steps in the negative direction, with the property 
that no two
of the lattice paths have a point in common. In the sequel, whenever
we use the term ``lattice path" we mean a lattice path 
consisting of horizontal unit steps in the positive direction 
and vertical unit steps in the negative direction. As is usually done, we
call a family of lattice paths with the property that no two
of the lattice paths have a point in common {\it nonintersecting\/}.

\midinsert

\vskip15pt
\vbox{
\centertexdraw{
\move (-1.5 0)
\bsegment
  \drawdim truecm  \linewd.08
  \move(0 0)
  \RhombusA \RhombusB \RhombusB 
  \RhombusA \RhombusA \RhombusB \RhombusA \RhombusA \RhombusA
  \RhombusB \RhombusB
  \move (-0.866025403784439 -.5)
  \RhombusA \RhombusB \RhombusB \RhombusB \RhombusB
  \RhombusA \RhombusA \RhombusB \RhombusA \RhombusA \RhombusB
  \move (-1.732050807568877 -1)
  \RhombusB \RhombusB \RhombusA \RhombusB \RhombusB \RhombusB
  \RhombusA \RhombusA \RhombusA \RhombusB \RhombusA
  \move (1.732050807568877 -4)
    \bsegment
  \rlvec(0.866025403784439 -.5) \rlvec(-0.866025403784439 -.5) 
  \lfill f:0
  \savepos(0 -1)(*ex *ey)
    \esegment
  \move(*ex *ey)
  \RhombusA \RhombusA \RhombusB \RhombusA \RhombusB
  \move (1.732050807568877 0)
  \RhombusC \RhombusC \RhombusC \RhombusC \RhombusC 
  \move (1.732050807568877 -1)
  \RhombusC \RhombusC \RhombusC \RhombusC \RhombusC 
  \move (3.464101615137755 -3)
  \RhombusC \RhombusC \RhombusC 
  \move (-0.866025403784439 -.5)
  \RhombusC
  \move (-0.866025403784439 -1.5)
  \RhombusC
  \move (0.866025403784439 -2.5)
  \RhombusC \RhombusC 
  \move (0.866025403784439 -3.5)
  \RhombusC 
  \move (0 -5)
  \RhombusC \RhombusC 
  \move (2.598076211353316 -7.5)
  \RhombusC 
  \move (-1.732050807568877 -3)
  \RhombusC 
  \move (-1.732050807568877 -4)
  \RhombusC 
  \move (-1.732050807568877 -5)
  \RhombusC 
  \move (-1.732050807568877 -6)
  \RhombusC \RhombusC \RhombusC \RhombusC
  \move (0.4330127018922194 .25)
  \linewd.05
  \lcir r:.1
  \move (-0.4330127018922194 -.25)
  \lcir r:.15
  \move (-1.299038105676658 -.75)
  \lcir r:.15
  \move (2.165063509461097 -4.75)
  \lcir r:.15
  \move (3.031088913245536 -9.25)
  \lcir r:.15
  \move (3.897114317029974 -8.75)
  \lcir r:.15
  \move (4.763139720814414 -8.25)
  \lcir r:.15
  \move (5.629165124598852 -7.75)
  \lcir r:.15
  \linewd.05
  \move (0.4330127018922194 .25)
  \hdSchritt \vdSchritt \vdSchritt \hdSchritt \hdSchritt \vdSchritt 
     \hdSchritt \hdSchritt \hdSchritt \vdSchritt \vdSchritt 
  \move (-0.4330127018922194 -.25)
  \hdSchritt \vdSchritt \vdSchritt \vdSchritt \vdSchritt \hdSchritt 
     \hdSchritt \vdSchritt \hdSchritt \hdSchritt \vdSchritt 
  \move (-1.299038105676658 -.75)
  \vdSchritt \vdSchritt \hdSchritt \vdSchritt \vdSchritt \vdSchritt 
     \hdSchritt \hdSchritt \hdSchritt \vdSchritt \hdSchritt
  \move (2.165063509461097 -4.75)
  \hdSchritt \hdSchritt \vdSchritt \hdSchritt \vdSchritt 
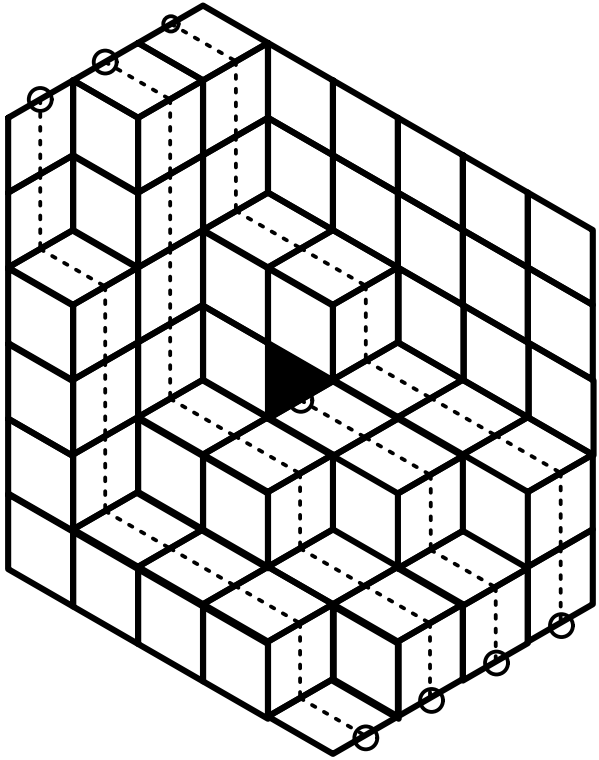
\htext (-1.3 -10.5){\eightpoint a. The paths associated to a rhombus tiling}
\esegment
\move (1.5 0)
\bsegment
  \drawdim truecm  \linewd.05
  \move (0.4330127018922194 .25)
  \lcir r:.15
  \move (-0.4330127018922194 -.25)
  \lcir r:.15
  \move (-1.299038105676658 -.75)
  \lcir r:.15
  \move (2.165063509461097 -4.75)
  \lcir r:.15
  \move (3.031088913245536 -9.25)
  \lcir r:.15
  \move (3.897114317029974 -8.75)
  \lcir r:.15
  \move (4.763139720814414 -8.25)
  \lcir r:.15
  \move (5.629165124598852 -7.75)
  \lcir r:.15
  \move (0.4330127018922194 .25)
  \rmove (0.866025403784439 -.5)
  \fcir f:0 r:.05
  \rmove (0 -1)
  \fcir f:0 r:.05
  \rmove (0 -1)
  \fcir f:0 r:.05
  \rmove (0.866025403784439 -.5)
  \fcir f:0 r:.05
  \rmove (0.866025403784439 -.5)
  \fcir f:0 r:.05
  \rmove (0 -1)
  \fcir f:0 r:.05
  \rmove (0.866025403784439 -.5)
  \fcir f:0 r:.05
  \rmove (0.866025403784439 -.5)
  \fcir f:0 r:.05
  \rmove (0.866025403784439 -.5)
  \fcir f:0 r:.05
  \rmove (0 -1)
  \fcir f:0 r:.05
  \move (-0.4330127018922194 -.25)
  \rmove (0.866025403784439 -.5)
  \fcir f:0 r:.05
  \rmove (0 -1)
  \fcir f:0 r:.05
  \rmove (0 -1)
  \fcir f:0 r:.05
  \rmove (0 -1)
  \fcir f:0 r:.05
  \rmove (0 -1)
  \fcir f:0 r:.05
  \rmove (0.866025403784439 -.5)
  \fcir f:0 r:.05
  \rmove (0.866025403784439 -.5)
  \fcir f:0 r:.05
  \rmove (0 -1)
  \fcir f:0 r:.05
  \rmove (0.866025403784439 -.5)
  \fcir f:0 r:.05
  \rmove (0.866025403784439 -.5)
  \fcir f:0 r:.05
  \move (-1.299038105676658 -.75)
  \rmove (0 -1)
  \fcir f:0 r:.05
  \rmove (0 -1)
  \fcir f:0 r:.05
  \rmove (0.866025403784439 -.5)
  \fcir f:0 r:.05
  \rmove (0 -1)
  \fcir f:0 r:.05
  \rmove (0 -1)
  \fcir f:0 r:.05
  \rmove (0 -1)
  \fcir f:0 r:.05
  \rmove (0.866025403784439 -.5)
  \fcir f:0 r:.05
  \rmove (0.866025403784439 -.5)
  \fcir f:0 r:.05
  \rmove (0.866025403784439 -.5)
  \fcir f:0 r:.05
  \rmove (0 -1)
  \fcir f:0 r:.05
  \move (2.165063509461097 -4.75)
  \rmove (0.866025403784439 -.5)
  \fcir f:0 r:.05
  \rmove (0.866025403784439 -.5)
  \fcir f:0 r:.05
  \rmove (0 -1)
  \fcir f:0 r:.05
  \rmove (0.866025403784439 -.5)
  \fcir f:0 r:.05
  \linewd.05
  \move (0.4330127018922194 .25)
  \hdSchritt \vdSchritt \vdSchritt \hdSchritt \hdSchritt \vdSchritt 
     \hdSchritt \hdSchritt \hdSchritt \vdSchritt \vdSchritt 
  \move (-0.4330127018922194 -.25)
  \hdSchritt \vdSchritt \vdSchritt \vdSchritt \vdSchritt \hdSchritt 
     \hdSchritt \vdSchritt \hdSchritt \hdSchritt \vdSchritt 
  \move (-1.299038105676658 -.75)
  \vdSchritt \vdSchritt \hdSchritt \vdSchritt \vdSchritt \vdSchritt 
     \hdSchritt \hdSchritt \hdSchritt \vdSchritt \hdSchritt
  \move (2.165063509461097 -4.75)
  \hdSchritt \hdSchritt \vdSchritt \hdSchritt \vdSchritt 
\htext (0 -10.5){\eightpoint b. The paths, isolated}
\esegment
}

\vskip.5cm

$$
\Einheit.64cm
\Gitter(9,9)(0,0)
\Koordinatenachsen(9,9)(0,0)
\Kreis(0,6)
\Kreis(1,7)
\Kreis(2,8)
\Kreis(4,4)
\Kreis(5,0)
\Kreis(6,1)
\Kreis(7,2)
\Kreis(8,3)
\Pfad(0,4),22\endPfad
\Pfad(0,4),1\endPfad
\Pfad(1,1),222\endPfad
\Pfad(1,1),111\endPfad
\Pfad(4,0),2\endPfad
\Pfad(4,0),1\endPfad
\Pfad(1,7),1\endPfad
\Pfad(2,3),2222\endPfad
\Pfad(2,3),11\endPfad
\Pfad(4,2),2\endPfad
\Pfad(4,2),11\endPfad
\Pfad(6,1),2\endPfad
\Pfad(2,8),1\endPfad
\Pfad(3,6),22\endPfad
\Pfad(3,6),11\endPfad
\Pfad(5,5),2\endPfad
\Pfad(5,5),111\endPfad
\Pfad(8,3),22\endPfad
\Pfad(4,4),11\endPfad
\Pfad(6,3),2\endPfad
\Pfad(6,3),1\endPfad
\Pfad(7,2),2\endPfad
\Label\l{P_1}(1,2)
\Label\l{P_2}(2,5)
\Label\o{P_3}(5,6)
\Label\u{P_4}(5,4)
\Label\lo{A_1}(0,6)
\Label\lo{A_2}(1,7)
\Label\lo{A_3}(2,8)
\Label\l{A_4}(4,4)
\Label\ru{E_1}(5,0)
\Label\ru{E_2}(6,1)
\Label\ru{E_3}(7,2)
\Label\ru{E_4}(8,3)
\hskip4.5cm
$$
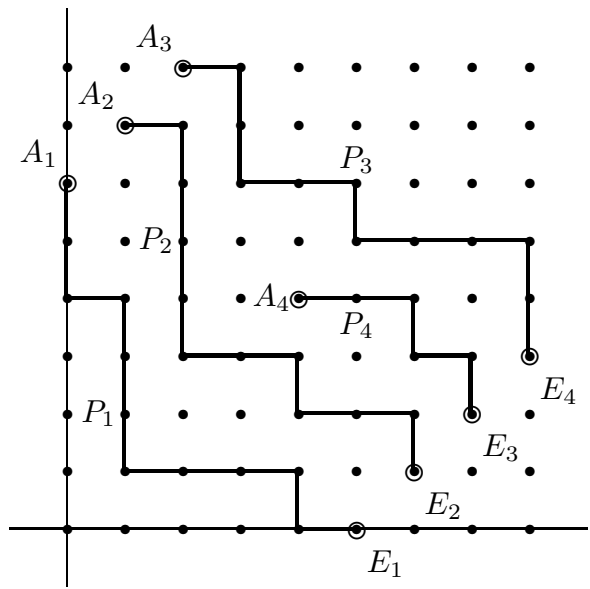
\centerline{\eightpoint c. After deformation: A family of
nonintersecting lattice paths}
\vskip10pt
\centerline{\eightpoint Figure 3}
}
\vskip10pt
\endinsert

It is easy to see that this correspondence sets up a bijection
between rhombus tilings of a hexagon with sides
$a,b+1,c,a+1,b,c+1$, of which the central triangle is removed,
and families $(P_1,P_2,\dots,P_{a+1})$ of nonintersecting lattice
paths, where for $i=1,2\dots,a$ the path $P_i$ runs from $A_i=(i-1,c+i)$ 
to one of the points $E_j=(b+j-1,j-1)$,
$j=1,2,\dots,a+1$, and where $P_{a+1}$ runs from
$A_{a+1}=\big((a+b)/2,(a+c)/2\big)$ to one of the points $E_j=(b+j-1,j-1)$,
$j=1,2,\dots,a+1$. It is the latter
enumeration problem that we are going to solve in the subsequent
sections. 

As mentioned in the Introduction, this enumeration problem
for nonintersecting lattice paths is unusual, as the starting points
are not lined up in the ``right" way (in this case this would mean
from bottom-left to top-right), the $(a+1)$-st starting point being
located in the middle of the region. So, it is not fixed which
starting point is connected to which end point. In our example of
Figure~3.c, $A_1$ is connected with $E_1$, $A_2$ is connected with
$E_2$, $A_3$ is connected with $E_4$, and the ``exceptional" starting
point $A_4$ is connected with $E_3$. But there are other
possibilities. In general, for any $k$, $1\le k\le a+1$, it is possible
that $A_1$ is connected with $E_1$, \dots, $A_{k}$ is connected with
$E_{k}$, $A_{k+1}$ is connected with $E_{k+2}$, \dots, $A_a$ is connected 
with $E_{a+1}$, and the ``exceptional" starting point $A_{a+1}$ is connected 
with $E_{k+1}$. To obtain the solution of our enumeration problem, we have to
find the number of all nonintersecting lattice paths in each case and
then form the sum over all $k$.

\subhead 3. From nonintersecting lattice paths to Schur functions\endsubhead
The aim of this section is to describe, starting from the interpretation of our
problem in
terms of nonintersecting lattice paths that was derived in the 
previous section, how to derive an expression for the number of
rhombus tilings that we are interested in. This expression, displayed
in (\CD), features (specialized) Schur functions.

Usually, when enumerating nonintersecting lattice paths, one obtains
a determinant by means of the main theorem on 
nonintersecting lattice paths \cite{\GeViAB, Cor.~2; \StemAE,
Theorem~1.2}. However, as mentioned in the Introduction, the
arrangement of the starting points in our problem, as described in
the previous section, does not allow a direct application of this theorem. 

\vskip10pt
\vbox{
$$
\Einheit.64cm
\Gitter(9,9)(0,0)
\Koordinatenachsen(9,9)(0,0)
\Kreis(0,6)
\Kreis(1,7)
\Kreis(2,8)
\Kreis(4,4)
\Kreis(5,0)
\Kreis(6,1)
\Kreis(7,2)
\Kreis(8,3)
\DickPunkt(0,0)
\DickPunkt(1,1)
\DickPunkt(2,2)
\DickPunkt(3,3)
\DickPunkt(5,5)
\DickPunkt(6,6)
\DickPunkt(7,7)
\DickPunkt(8,8)
\Pfad(0,4),22\endPfad
\Pfad(0,4),1\endPfad
\Pfad(1,1),222\endPfad
\Pfad(1,1),111\endPfad
\Pfad(4,0),2\endPfad
\Pfad(4,0),1\endPfad
\Pfad(1,7),1\endPfad
\Pfad(2,3),2222\endPfad
\Pfad(2,3),11\endPfad
\Pfad(4,2),2\endPfad
\Pfad(4,2),11\endPfad
\Pfad(6,1),2\endPfad
\Pfad(2,8),1\endPfad
\Pfad(3,6),22\endPfad
\Pfad(3,6),11\endPfad
\Pfad(5,5),2\endPfad
\Pfad(5,5),111\endPfad
\Pfad(8,3),22\endPfad
\Pfad(4,4),11\endPfad
\Pfad(6,3),2\endPfad
\Pfad(6,3),1\endPfad
\Pfad(7,2),2\endPfad
\Label\l{P_1}(1,2)
\Label\l{P_2}(2,5)
\Label\o{P_3}(5,6)
\Label\u{P_4}(5,4)
\hskip4.5cm
$$
\vskip10pt
\centerline{\eightpoint Figure 4}
}
\vskip10pt

In order to access the problem, as formulated at the end of the
previous section, we set it up as follows. Each of the lattice paths $P_i$,
$i=1,2,\dots,a$, has to pass through a lattice point on the
diagonal line $x-y=(b-c)/2$, the central point
$A_{a+1}=\big((a+b)/2,(a+c)/2\big)$ (which is the starting point for
path $P_{a+1}$) excluded (see Figure~4; the points on
the diagonal where the paths may pass through are 
indicated by bold dots). So, as in the previous section, let $k$ be
an integer between $1$ and $a+1$. For any fixed choice of points
$M_1,\dots,M_{k},M_{k+2},\dots,M_{a+1}$ on the diagonal $x-y=(b-c)/2$, 
such that $M_1$ is to the
left of $M_2$, \dots, $M_{k}$ is to the left of the central point
$A_{a+1}$, $A_{a+1}$ is to the left of $M_{k+2}$, \dots, $M_{a}$ is to
the left of $M_{a+1}$, the number of 
families $(P_1,P_2,\dots,P_{a+1})$ of nonintersecting lattice
paths, where for $i=1,2\dots,k$ the path $P_i$ runs from $A_i$
through $M_i$ to $E_i$, 
where for $i=k+1,k+2,\dots,a$ the path $P_i$ runs from $A_i$
through $M_{i+1}$ to $E_{i+1}$, and where $P_{a+1}$ runs from
$A_{a+1}$ to $E_{k+1}$, is easily computed using the main theorem on
nonintersecting lattice paths. For convenience, let $M_{k+1}=A_{a+1}$.
Then this number equals
$$\det_{1\le i,j\le a}\(\vert \Cal P(A_i\to M_{j+\chi(j\ge k+1)})\vert\)
\cdot\det_{1\le i,j\le a+1}\(\vert \Cal P(M_i\to E_{j})\vert\),
\tag\CA$$
where $\vert \Cal P(A\to B)\vert$ denotes the number of all lattice
paths from $A$ to $B$, and where $\chi$ is the usual truth function,
$\chi(\Cal A)=1$ if $\Cal A$ is true and $\chi(\Cal A)=0$ otherwise.

It is now no difficulty to observe that the expression (\CA) can be
rewritten using (specialized) Schur functions. For information on
Schur functions and related definitions 
we refer the reader to Chapter~I of Macdonald's
classical book \cite{\MacdAC}. There are many ways to express a Schur
function. The one we need here is the N\"agelsbach--Kostka formula 
(the ``dual'' Jacobi-Trudi identity; see \cite{\MacdAC, Ch.~I, (3.5)}). Let
$\la$ be a partition with largest part at most $m$. Then the Schur
function $s_\la(x_1,x_2,\dots,x_n)$ is given by
$$s_\la(x_1,x_2,\dots,x_n)=\det_{1\le i,j\le m}\big(e_{{}(\trans\la)_{i}-i+j}
(x_1,x_2,\dots,x_n)\big),$$
where $e_s(x_1,\dots,x_n):=\sum _{1\le i_1<\dots< i_s\le n}
^{}x_{i_1}\cdots x_{i_s}$ is the {\it elementary symmetric
function\/} of order $s$ in $x_1,\dots,x_n$, and
where $\trans\la$ denotes the partition conjugate to $\la$. We 
write briefly $X_n$ for the set of variables 
$\{x_1,x_2,\dots,x_n\}$. In particular, the symbol $s_\la(X_n)$ will
be short for $s_\la(x_1,x_2,\dots,x_n)$.

Now suppose that $M_\ell=\big((a+b)/2+i_\ell,(a+c)/2+i_\ell\big)$,
$\ell=1,2,\dots,a+1$. In particular, since $M_{k+1}=A_{a+1}=
\big((a+b)/2,(a+c)/2\big)$, we have $i_{k+1}=0$. The ordering of the
points $M_\ell$ from left to right as $\ell$ increases implies 
$-(a+b)/2 \le i_1 < \cdots < i_k < i_{k+1}=0
 < i_{k+2} < \cdots < i_{a+1} \le (a+b)/2$.
Then it is straight-forward to see that, using the above notation, 
the expression (\CA) equals
$$s_\la(X_{(b+c+2)/2})\cdot s_\mu(X_{(b+c)/2})
\Big\vert_{x_1=x_2=\dots=x_{(b+c+2)/2}=1}\ ,\tag\CB$$
where $\la$ and $\mu$ are the partitions whose conjugates $\trans\la$ and
$\trans\mu$ are given by
$$
\matrix\format\r&\l&\quad \quad \l\\
\trans \lambda_h &=
 \dfrac{b-a}{2} +i_{a+1-h+\chi(a+1-h\ge k+1)} +h
 &\text{for $h=1, \dots, a$,} \\
\trans \mu_h &=
 \dfrac{b-a}{2} - i_h +h-1
 &\text{for $h=1, \dots, a+1$.}
\endmatrix
\tag\CC
$$
(In view of the geometric interpretation of the N\"agelsbach--Kostka
formula in terms of nonintersecting lattice paths, see \cite{\FuKrAA,
Sec.~4, Fig.~8; \SagaAL, Sec.~4.5},
this could also be read off directly from the lattice path picture of
our setup, as exemplified by Figure~4.)

According to the preceding considerations, what we would like to do
in order to find the {\it total\/} number of 
all rhombus tilings of a hexagon with sides
$a,b+1,c,a+1,b,c+1$, of which the central triangle is removed, is to
sum the products (\CB) of Schur functions over all possible choices
of $\la$ and $\mu$. That is, we want to find the sum
$$\sum_{(\la,\mu)} \big(s_\la(X_{(b+c+2)/2})\cdot
s_\mu(X_{(b+c)/2})\big)\Big\vert_{x_1=x_2=\dots=x_{(b+c+2)/2}=1}\ ,$$
where $(\la,\mu)$ ranges over all possible pairs of partitions such
that (\CC) is satisfied, for some $k$ and $i_1,i_2,\dots,i_{a+1}$ as
above. 

Summarizing, we have shown the following.

\proclaim{Proposition~\TB}Let $a,b,c$ be positive integers, all of the
same parity. Then the number of all rhombus tilings of a hexagon with sides
$a,b+1,c,a+1,b,c+1$, of which the central triangle is removed, equals
$$\sum_{(\la,\mu)} \big(s_\la(X_{(b+c+2)/2})\cdot
s_\mu(X_{(b+c)/2})\big)\Big\vert_{x_1=x_2=\dots=x_{(b+c+2)/2}=1}\ ,\tag\CD$$
where $(\la,\mu)$ ranges over all possible pairs of partitions such
that {\rm(\CC)} is satisfied, for some $k$ with $1\le k\le a+1$, and for
some $i_1,i_2,\dots,i_{a+1}$ with 
$-(a+b)/2 \le i_1 < \cdots < i_k < i_{k+1}=0
 < i_{k+2} < \cdots < i_{a+1} \le (a+b)/2$.
\endproclaim

\subhead 4. The main theorem and its implications\endsubhead
As a matter of fact, even the {\it
unspecialized\/} sum that appears in (\CD) can be
evaluated in closed form. This is the subject of the subsequent
Theorem~\TC, which is the main theorem of our paper. 
Theorem~\TA\ of the Introduction then follows
immediately. Before we proceed to the proof of Theorem~\TC\ in
Section~5, in this section 
we formulate a conjectured generalization of Theorem~\TC\ as
Conjecture~\TE,
and we complement the enumeration result in Theorem~\TA\ by a further
result (Theorem~\TD) on the enumeration of rhombus tilings
of a ``punctured" hexagon, which follows also from Theorem~\TC.

\proclaim{Theorem~\TC}Let $a,b,n$ be positive integers, where $a$ and
$b$ are of the same parity. Then
$$\multline
\sum_{(\la,\mu)} s_\la(X_{n+1})\cdot
s_\mu(X_{n})\\
=s_{(\cl{(a+1)/2}^{\fl{b/2}})}(X_n)\cdot
s_{(\cl{a/2}^{\cl{b/2}})}(X_n)\cdot
s_{(\cl{a/2}^{\cl{b/2}})}(X_{n+1})\cdot
s_{(\fl{a/2}^{\cl{(b+1)/2}})}(X_{n+1}),
\endmultline\tag\DA$$
where $(\la,\mu)$ ranges over all possible pairs of partitions such
that {\rm(\CC)} is satisfied, for some $k$ with $1\le k\le a+1$, and for
some $i_1,i_2,\dots,i_{a+1}$ with 
$-(a+b)/2 \le i_1 < \cdots < i_k < i_{k+1}=0
 < i_{k+2} < \cdots < i_{a+1} \le (a+b)/2$.
\endproclaim

Theorem~\TA\ follows immediately from Proposition~\TB\ and 
the $n=(b+c)/2$ special case of
Theorem~\TC\ by using the well-known fact (see \cite{\MacdAC, Ch.~I,
Sec.~5, Ex.~13.(b), $q\to1$}) 
that the evaluation of a
Schur function $s_{(A^B)}(X_n)$ of rectangular shape at
$x_1=x_2=\dots=x_{n}=1$ counts the number of all plane partitions
inside the $A\times B\times (n-B)$ box.

A different specialization of Theorem~\TC\ leads to an
enumeration result for rhombus tilings of a ``punctured"
hexagon with sides $a,b+1,c,a+1,b,c+1$ where {\it not all\/} of $a,b,c$ are
of the same parity. Namely, on setting $x_1=x_2=\dots=x_n=1$, 
$x_{n+1}=0$, $n=(b+c)/2$ in Theorem~\TC, and on finally replacing
$c$ by $c+1$, we obtain the following result by performing
the analogous translations that lead from
rhombus tilings to nonintersecting lattice paths and finally to
specialized Schur functions.

\vskip15pt
\vbox{
\centertexdraw{
  \drawdim truecm  \linewd.02
  \move (-0.866025403784439 -.5)
  \rmove (0.866025403784439 -.5)
  \rhombus \rhombus 
  \rhombus \rhombus \rhombus \rhombus \ldreieck
  \move (-1.732050807568877 -1)
  \rmove (0.866025403784439 -.5)
  \rhombus \rhombus \rhombus 
  \rhombus \rhombus \rhombus \rhombus \ldreieck
  \move (-1.732050807568877 -2)
  \rhombus \rhombus \rhombus \rhombus 
  \rhombus \rhombus \rhombus \rhombus \ldreieck
  \move (-1.732050807568877 -2)
  \rdreieck
  \rhombus \rhombus \rhombus \rhombus 
  \rhombus \rhombus \rhombus \rhombus \ldreieck
  \move (-1.732050807568877 -3)
  \rdreieck
  \rhombus \rhombus \rhombus \rhombus 
  \rhombus \rhombus \rhombus \rhombus 
  \move (-1.732050807568877 -4)
  \rdreieck
  \rhombus \rhombus \rhombus \rhombus 
  \rhombus \rhombus \rhombus 
  \move (-1.732050807568877 -5)
  \rdreieck
  \rhombus \rhombus \rhombus \rhombus 
  \rhombus \rhombus 
  \move (-1.732050807568877 -6)
  \rdreieck
  \rhombus \rhombus \rhombus \rhombus 
  \rhombus 
  \move (1.732050807568877 -4)
    \bsegment
  \rlvec(0.866025403784439 -.5) \rlvec(-0.866025403784439 -.5) 
  \lfill f:0
  \savepos(0 -1)(*ex *ey)
    \esegment
  \move(*ex *ey)
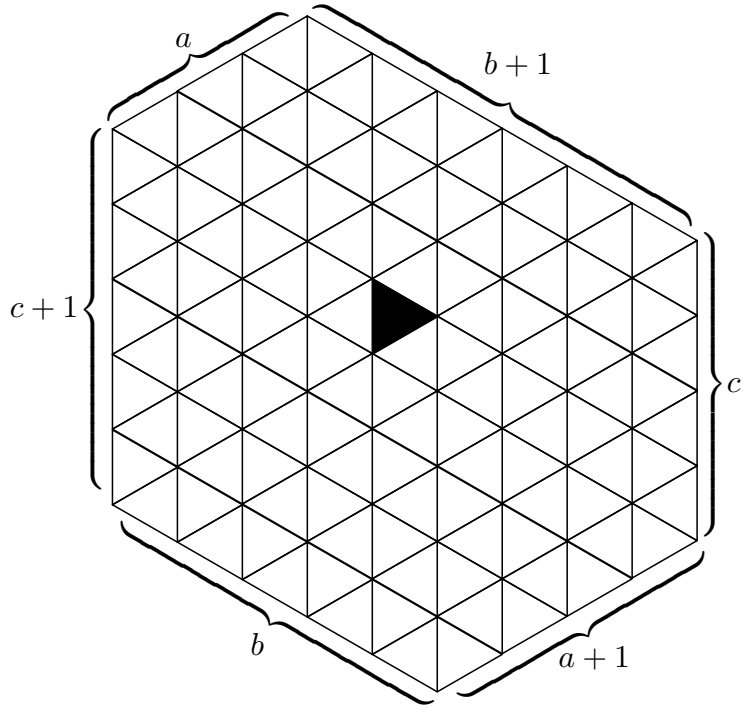
\htext (-2.8 -10.5){\eightpoint A ``punctured" hexagon with sides $a,b+1,c,a+1,b,c+1$,}
\htext (-2.3 -11){\eightpoint $c$ being of parity different from the parity of $a$ and $b$}
\rtext td:0 (6 -5.5){$\sideset {} \and c\to 
    {\left.\vbox{\vskip1.8cm}\right\}}$}
\rtext td:60 (3.46 -2){$\sideset {} \and {} \to 
    {\left.\vbox{\vskip2.6cm}\right\}}$}
\rtext td:120 (-.44 -1.25){$\sideset {} \and {} \to 
    {\left.\vbox{\vskip1.3cm}\right\}}$}
\rtext td:0 (-3.1 -4.5){$\sideset {c+1} \and {}\to 
    {\left\{\vbox{\vskip2.1cm}\right.}$}
\rtext td:240 (0.44 -8.25){$\sideset {} \and {} \to 
    {\left.\vbox{\vskip2.1cm}\right\}}$}
\rtext td:300 (4.33 -8.5){$\sideset {} \and {} \to 
    {\left.\vbox{\vskip1.7cm}\right\}}$}
\htext (-.9 -.9){$a$}
\htext (3.2 -1.3){$b+1$}
\htext (4.2 -9.2){$a+1$}
\htext (0.1 -9){$b$}
}
\vskip4pt
\centerline{\eightpoint Figure 5}
}
\vskip10pt

\proclaim{Theorem~\TD}Let $a,b,c$ be positive integers, $a$ and $b$
of the same parity, $c$ of different parity. From the hexagon 
with sides $a,b+1,c,a+1,b,c+1$ the triangle is removed, 
whose side parallel to the sides of lengths $a$ and $a+1$ lies on the
line that is equidistant to these sides, 
whose side parallel to the sides of lengths $b$ and $b+1$ lies on the
line that is by two ``units" closer to the side of length $b+1$ than
to the side of length $b$, and
whose side parallel to the sides of lengths $c$ and $c+1$ lies on the
line that is by one ``unit" closer to the side of length $c+1$ than
to the side of length $c$, 
see Figure~5.
Then the number of all rhombus tilings of this ``punctured" 
hexagon with sides $a,b+1,c,a+1,b,c+1$ equals
$$
B\(\fl{\frac {\vphantom{b}a+2} {2}},\fl{\frac {b} {2}},\fl{\frac
{\vphantom{b}c+2} {2}}\)
B\(\fl{\frac {a+1} {2}},\fl{\frac {b+1} {2}},\fl{\frac {\vphantom{b}c+1}
{2}}\)^2
B\(\fl{\frac {\vphantom{b}a} {2}},\fl{\frac {b+2} {2}},\fl{\frac
{\vphantom{b}c}
{2}}\),
\tag\DB$$
where $B(\alpha,\beta,\gamma)$ is the number of all plane partitions inside the
$\alpha\times \beta\times \gamma$ box, which is given by \rm(\AA).
\endproclaim

Even if Theorem~\TC\ suffices to prove Theorem~\TA\ (and
Theorem~\TD), it is certainly of interest that apparently there 
holds a symmetric generalization of Theorem~\TC.
We have overwhelming evidence through computer computations that
actually the following is true.
\proclaim{Conjecture~\TE}Let $a,b,n$ be positive integers, where $a$ and
$b$ are of the same parity. Then
$$\multline
\sum_{(\la,\mu)} s_\la(X_{n+2})\cdot
s_\mu(X_{n})
=s_{(\cl{(a+1)/2}^{\fl{b/2}})}(X_n)\cdot
s_{(\cl{a/2}^{\cl{b/2}})}(X_n\cup\{x_{n+1}\})\\
\times
s_{(\cl{a/2}^{\cl{b/2}})}(X_{n}\cup\{x_{n+2}\})\cdot
s_{(\fl{a/2}^{\cl{(b+1)/2}})}(X_{n+2}),
\endmultline\tag\DC$$
where $(\la,\mu)$ ranges over all possible pairs of partitions such
that {\rm(\CC)} is satisfied, for some $k$ with $1\le k\le a+1$, and for
some $i_1,i_2,\dots,i_{a+1}$ with 
$-(a+b)/2 \le i_1 < \cdots < i_k < i_{k+1}=0
 < i_{k+2} < \cdots < i_{a+1} \le (a+b)/2$.
\endproclaim

We remark that the specialization $n=(b+c)/2$,
$x_1=x_2=\dots=x_{(b+c+4)/2}=1$ of Conjecture~\TE\ leads, 
up to parameter permutation,
again to Theorem~\TD, thus providing further evidence for the truth
of the Conjecture.

\subhead 5. The minor summation formula and proof of the main
theorem\endsubhead
This final section is devoted to a proof of Theorem~\TC. It makes essential
use of the minor summation formula of Ishikawa and Wakayama. An
outline of the proof is as follows. First, the minor summation
formula is used to convert the sum on the left-hand side of (\DA)
into a Pfaffian. This Pfaffian can be easily reduced to a
determinant. In Lemma~\TI\ it is seen that this determinant factors
into a product of two Pfaffians. Finally, Lemma~\TJ\ shows that each of
these Pfaffians is, basically, a product of two Schur functions, so
that the final form of the result, as given by the right-hand side
of (\DA), follows from a simple computation, by which we conclude this
section.

To begin with, let us recall the minor summation formula due to Ishikawa
 and Wakayama \cite{\IsWaAA, Theorem~2}.

\proclaim{Theorem~\TG}
Let $n$, $p$, $q$ be integers such that $n+q$ is even and $0 \le n-q \le p$.
Let $G$ be any $n \times p$ matrix, $H$ be any $n \times q$ matrix,
 and $A=(a_{ij})_{1 \le i,j \le p}$ be any skew-symmetric matrix.
Then we have
$$
\sum_{K}
 \Pf \left( A^K_K \right)
 \det \big( G_K \ \vdots \  H \big)
=
(-1)^{q(q-1)/2}
 \Pf \left( \matrix
 G\,A\,\,\trans G & H \\ - \trans H & 0 \endmatrix \right),
\tag\EE$$
where $K$ runs over all $(n-q)$-element subsets of $[1, p]$,
 $A^K_K$ is the skew-symmetric matrix obtained by picking the rows and
 columns indexed by $K$ and $G_K$ is the sub-matrix of $G$ consisting of
 the columns corresponding to $K$.
\endproclaim

In order to apply this formula, we have to first describe the pairs
$(\la,\mu)$ of partitions over which the sum in (\DA) is taken
{\it directly}. This is done in Lemma~\TH\ below. 
The reader should be reminded that 
the description that the formulation of Theorem~\TC\
gives is in terms of the {\it conjugates} of $\la$ and $\mu$ (compare
(\CC)), which is not suitable for application of the minor
summation formula (\EE).

For convenience,
let $\Cal{R}(a,b)$ be the set of all pairs $(\lambda, \mu)$ of partitions 
satisfying (\CC), for some $k$ with $1\le k\le a+1$, and for
some $i_1,i_2,\dots,i_{a+1}$ with 
$-(a+b)/2 \le i_1 < \cdots < i_k < i_{k+1}=0
 < i_{k+2} < \cdots < i_{a+1} \le (a+b)/2$.

First note that $l(\lambda) \le b+1$ and $l(\mu) \le b$ for $(\lambda, \mu)
 \in \Cal{R}(a,b)$.

\proclaim{Lemma~\TH}
If we define subsets $J = \{ j_1 < \cdots < j_{b+1} \}$ and $J'
 = \{ j'_1 < \cdots < j'_b \}$ by the relations
$$
\align
j_h &= \lambda_{b+2-h} + h-1 \quad\text{for $h=1, \dots, b+1$,} \\
j'_h &= \mu_{b+1-h} + h-1 \quad\text{for $h=1, \dots, b$,}
\endalign
$$
then we have
\roster
\item
$(a+b)/2 \in J$.
\item
$J' = \{ a+b - j : j \in J, j \neq (a+b)/2 \}$.
\endroster
\endproclaim

\demo{Proof}
(1)
Since $i_{k+2} \ge 1$ and $i_k \le -1$, we have
$$
\align
{}^t \lambda_{a-k} &= \frac{b-a}{2}+i_{k+2}+a-k \ge \frac{a+b}{2} - k+1, \\
{}^t \lambda_{a-k+1} &= \frac{b-a}{2}+i_k+a-k+1 \le \frac{a+b}{2} - k.
\endalign
$$
Hence we obtain $\lambda_{(a+b)/2-k+1} = a-k$, which means
 $j_{(b-a)/2+k+1} = (a+b)/2$.

(2)
It is enough to show that
$$
\align
&\lambda_{b+2-h} + \mu_h = a+1
 \quad\text{for $h=1, \dots, \frac{b-a}{2}+k$,} \\
&\lambda_{b+2-h} + \mu_{h-1} = a
 \quad\text{for $h= \frac{b-a}{2}+k+2, \dots, b+1$.}
\endalign
$$
By (\CC),
the parts of $\lambda$ and $\mu$ can be expressed in terms of the $i_j$'s:
$$
\align
\lambda_{b+2-h} &=
 \# \left\{ j : i_{a+1-j+\chi(a+1-j\ge k+1)} \ge \frac{a+b}{2}+2-h-j \right\}, \\
\mu_h &=
 \# \left\{ j : i_j \le \frac{b-a}{2}-1-h+j \right\}.
\endalign
$$
If $h \le (b-a)/2+k$, then $i_k < (a+b)/2+2-h-(a+1-k)$ and so we have
 $\lambda_{b+2-h} \le a-k$.
This implies
$$
\lambda_{b+2-h} = \# \left\{ j : i_j \ge \frac{b-a}{2}-h+j \right\}
$$
and $\lambda_{b+2-h} + \mu_h = a+1$.
Similarly, since $\mu_{h-1} \le k$ for $h \ge (b-a)/2+k+2$, we have
$$
\align
\lambda_{b+2-h} &= \# \left\{ j : i_{j+\chi(j\ge k+1)} \ge \frac{b-a}{2}-h+j+1 \right\}, \\
\mu_{h-1} &= \# \left\{ j : i_{j+\chi(j\ge k+1)} \le \frac{b-a}{2}-h+j \right\},
\endalign
$$
and $\lambda_{b+2-h}+ \mu_{h-1} = a$.\quad \quad \qed
\enddemo

Now we describe our choices of matrices $G$, $H$ and $A$ for our
application of the minor summation formula in Theorem~\TG.
For a subset $I = \{ i_1 < \cdots < i_p \}$ of nonnegative integers,
 let $M_I(X_n)$ be the $n \times p$ matrix with $(k,l)$ entry $x_k^{i_l}$.
We define three subsets $P$, $Q$, $R$ as follows:
$$
\align
P &=
 \left\{ n-b, n-b+1, \dots, n-b+\frac{a+b}{2}-1, n-b+\frac{a+b}{2}+1, \dots,
 n-b+(a+b) \right\}, \\
Q &=
 \left\{ 0, 1, \dots, n-b-1, n-b+\frac{a+b}{2} \right\}, \\
R &=
 \{ 0, 1, \dots, n-b-1 \}.
\endalign
$$
Then let the matrices $G$ and $H$ be given by
$$
\align
G &= \pmatrix M_P(X_{n+1}) & 0 \\
              0            & M_P(X_n) \endpmatrix, \\
H &= \pmatrix M_Q(X_{n+1}) & 0 \\
              0            & M_R(X_n) \endpmatrix.
\endalign
$$
Furthermore, we define
$$
\Gamma =
 \left\{ 0, 1, \dots, \frac{a+b}{2}-1, \frac{a+b}{2}+1, \dots, a+b
\right\}.
$$
Let $A$ be the skew-symmetric matrix whose rows and columns are indexed by
 the set
$$
\Gamma \cup \bar{\Gamma} = \Gamma \cup \{ \bar{k} : k \in \Gamma \}
$$
and whose nonzero entries are given by
$$
a_{k,\bar{a+b-k}} = \cases
 \hphantom{-}1 &\text{if $1 \le k \le (a+b)/2-1$} \\
 -1 &\text{if $(a+b)/2+1 \le k \le a+b$.}
\endcases
$$
We apply Theorem~\TG\ to these matrices $G$, $H$ and $A$.

For a $2b$-element subset $K$ of $\Gamma \cup \bar{\Gamma}$, the sub-Pfaffian
 $\Pf(A^K_K)$ is easily computed.
This sub-Pfaffian vanishes unless the number of the unbarred elements in $K$
 is equal to that of the barred elements in $K$.
If $K = \{ j_1, \dots, j_b, \bar{j'_1}, \dots, \bar{j'_b} \}$, then
 we have
$$
\Pf (A^K_K) = \cases
 (-1)^{\# \{ h : j_h \ge (a+b)/2+1 \}} &\text{if $j_h+ j'_{b+1-h} = a+b$} \\
 0 &\text{otherwise.}
\endcases
$$
Now recall the bideterminantal
expression for Schur functions (see \cite{\MacdAC, Ch.~I, (3.1)}),
$$
s_\la(X_n)=s_\lambda (x_1, x_2,\dots ,x_n)=\frac { 
 \det\limits _{1\le i,j\le n}(x_j^{\lambda _i+n-i})} { 
\det\limits_{1\le i,j\le n}(x_j^{n-i})}=
\frac { 
 \det\limits _{1\le i,j\le n}(x_j^{\lambda _i+n-i})} {\Delta(X_n)}.\tag \EF
$$
Then for a subset $K$ as above, we have
$$
\align
\det ( G_K | H ) 
& = (-1)^{2b(n-b)+b+\# \{ h: j_h \ge (a+b)/2+1 \} +b(n-b)}
   s_{\lambda}(X_{n+1}) s_{\mu}(X_n)
   \Delta(X_{n+1}) \Delta(X_n) \\
& = (-1)^{bn+\# \{ h : j_h \ge (a+b)/2+1 \} }
   s_{\lambda}(X_{n+1}) s_{\mu}(X_n)
   \Delta(X_{n+1}) \Delta(X_n), \\
\endalign
$$
where $\lambda$ (resp. $\mu$) is the partition corresponding to the subset
 $J = \{ j_1, \dots, j_b, (a+b)/2 \}$ (resp. $J' = \{ a+b-j_1, \dots,
 a+b-j_b \}$) and $\Delta(X_n) = \prod_{1 \le i < j \le n} (x_j -x_i)$.
Therefore, if we apply the minor summation formula in Theorem~\TG, 
we obtain
$$
\align
\sum_{(\lambda, \mu) \in \Cal{R}(a,b)} s_{\lambda}(X_{n+1}) s_{\mu}(X_n)
&=
 \frac{(-1)^{bn}}{\Delta(X_{n+1}) \Delta(X_n)}
  \sum_{K \subset \Gamma \cup \bar{\Gamma}, |K|=2b}
  \Pf(A^K_K) \det (G_k | H) \\
&=
 \frac{(-1)^{bn+(2n-2b+1)(2n-2b)/2}}
      {\Delta(X_{n+1}) \Delta(X_n)}
  \Pf \left( \matrix
  G\,A\,\,\trans G & H \\ - \trans H & 0 \endmatrix \right) \\
&=
 \frac{(-1)^{bn+n-b}}
      {\Delta(X_{n+1}) \Delta(X_n)}
  \Pf \left( \matrix
  G\,A\,\,\trans G & H \\ - \trans H & 0 \endmatrix \right).
\tag{\EA}
\endalign
$$

If we write $A = \pmatrix 0 & B \\ B & 0\endpmatrix$, then we have
$$
\align
&\pmatrix G \,A \,\,\trans G & H \\
         - \trans H   & 0 \endpmatrix \\
&=
\pmatrix
 0       & M_P(X_{n+1}) \,B \,\,\trans M_P(X_n) & M_Q(X_{n+1}) & 0 \\
 -M_P(X_n) \,B \,\,\trans M_P(X_{n+1}) & 0      & 0            & M_R(X_n) \\
 - \trans M_Q(X_{n+1}) & 0                & 0            & 0 \\
 0       & - \trans M_R(X_n)              & 0            & 0
\endpmatrix.
\endalign
$$
By permuting rows and columns by the permutation
$$
\left( \matrix
1 & \cdots & n+1 & n+2 & \cdots & 3n-b+2 & 3n-b+3 & \cdots & 4n-2b+2 \\
1 & \cdots & n+1 & 2n-b+2 & \cdots & 4n-2b+2 & n+2 & \cdots & 2n-b+1
 \endmatrix \right),
$$
we have
$$
\align
&\Pf
\pmatrix G \,A \,\,\trans G & H \\
         - \trans H   & 0 \endpmatrix
\\
&=
(-1)^{(2n-b+1)(n-b)} \\
&\qquad \times
\Pf
\pmatrix
 0 & 0 & M_P(X_{n+1}) \,B \,\,\trans M_P(X_n) & M_Q(X_{n+1}) \\
 0 & 0 & - \trans M_R(X_n) & 0 \\
 -M_P(X_n) \,B \,\,\trans M_P(X_{n+1}) & M_R(X_n) & 0 & 0 \\
 - \trans M_Q(X_{n+1}) & 0 & 0 & 0
\endpmatrix
\\
&=
(-1)^{bn-n+(2n-b+1)(2n-b)/2}
 \det
\pmatrix
M_P(X_{n+1}) \,B \,\,\trans M_P(X_n) & M_Q(X_{n+1}) \\
- \trans M_R(X_n) & 0
\endpmatrix
\\
&=
(-1)^{bn+b(b-1)/2+(n-b)}
\det
\pmatrix
N(X_n, X_n) & M_R(X_n) & M_{\{ n-b+(a+b)/2 \}}(X_n) \\
- \trans M_R(X_n) & 0 & 0 \\
N(x_{n+1}, X_n) & M_R(x_{n+1}) & x_{n+1}^{n-b+(a+b)/2}
\endpmatrix
 \tag{\EB}
\endalign
$$
where
$$
N(X_m, Y_n) = M_P(X_m) \,B \,\,\trans M_P(Y_n).
$$
By direct computation, the $(i,j)$-entry of $N(X_m, Y_n)$ is equal to
$$
\frac
{x_i^{n-b} y_j^{n-b} (y_j^{(a+b)/2} - x_i^{(a+b)/2})
                     (y_j^{(a+b)/2+1} - x_i^{(a+b)/2+1})}
{y_j - x_i}.
$$

Here, the last determinant in \thetag{\EB} can be decomposed into the product of
 two Pfaffians by using the following lemma.

\proclaim{Lemma~\TI}
Let $A$ be an $n \times n$ skew-symmetric matrix,
 $b={}^t (b_1, \dots, b_n)$ and $c={}^t (c_1, \dots, c_n)$
 be column vectors, and $d$ a scalar.
Then the determinant of the $(n+1) \times (n+1)$ matrix
$$
\tilde{A} = \pmatrix
 A       & b \\
 -{}^t c & d \endpmatrix
$$
decomposes into the product of two Pfaffians as follows:
\roster
\item
If $n$ is even, then
$$
\det \tilde{A} =
 -
 \Pf (A)
 \Pf \pmatrix A       & b & c \\
              -{}^t b & 0 & -d \\
              -{}^t c & d & 0 \endpmatrix.
\tag\EG
$$
\item
If $n$ is odd, then
$$
\det \tilde{A} =
 \Pf \pmatrix A       & b  \\
              -{}^t b & 0  \endpmatrix
 \Pf \pmatrix A       & c  \\
              -{}^t c & 0  \endpmatrix.
\tag\EH
$$
\endroster
\endproclaim

\demo{Proof}
Expanding along the last column and the bottom row, we see that
$$
\align
\det \tilde{A}
&=
\sum_{i,j =1}^n (-1)^{n+1+i-2} b_i \cdot (-1)^{n+j-2} (-c_j)
  \cdot \det \left( A_{\hat{i}}^{\hat{j}} \right) 
+ d \det (A) \\
&=
\sum_{i,j =1}^n (-1)^{i+j} b_i c_j
   \det \left( A_{\hat{i}}^{\hat{j}} \right) 
+ d \det (A),
\endalign
$$
where $A_{\hat{i}}^{\hat{j}}$ denotes the matrix obtained from $A$
 by deleting the $i$-th row and the $j$-th column.

First suppose that $n$ is even.
Then, since $A_{\hat{i}}^{\hat{i}}$ (resp. $A$) is a skew-symmetric matrix
 of odd (resp. even) degree, $\det \left( A_{\hat{i}}^{\hat{i}} \right) =0$
 (resp. $\det(A) = \Pf(A)^2$).
By induction hypothesis, we see that, if $i < j$, then
$$
\align
\det \left( A_{\hat{i}}^{\hat{j}} \right) 
&=
(-1)^{(n-i-1)+(n-j)} \det
 \pmatrix
       &      &     & a_{1i} \\
       & A_{\hat{i},\hat{j}}^{\hat{i},\hat{j}} && \vdots \\
       &      &     & a_{ni} \\
a_{j1} & \cdots & a_{jn} & a_{ji}
\endpmatrix
\\
&=
(-1)^{i+j-1} \cdot (-1)
\Pf \left( A_{\hat{i},\hat{j}}^{\hat{i}, \hat{j}} \right)
\Pf
 \pmatrix
       &      &     & a_{1i} & a_{1j} \\
       & A_{\hat{i},\hat{j}}^{\hat{i},\hat{j}} && \vdots & \vdots \\
       &      &     & a_{ni} & a_{nj} \\
a_{i1} & \cdots & a_{in} & 0      & a_{ij} \\
a_{j1} & \cdots & a_{jn} & a_{ji} & 0 
\endpmatrix
\\
&=
(-1)^{i+j} \cdot (-1)^{(n-i-1)+(n-j)} 
\Pf \left( A_{\hat{i},\hat{j}}^{\hat{i}, \hat{j}} \right)
\Pf (A)
\\
&=
-
\Pf \left( A_{\hat{i},\hat{j}}^{\hat{i}, \hat{j}} \right)
\Pf (A).
\endalign
$$
Similarly, if $i>j$, then we have
$$
\det \left( A_{\hat{i}}^{\hat{j}} \right)
=
\Pf \left( A_{\hat{i},\hat{j}}^{\hat{i}, \hat{j}} \right)
\Pf (A).
$$
Hence we have
$$
\align
&\det \tilde{A} \\
&=
 \Pf(A)
 \left(
  \sum_{1 \le i < j \le n} (-1)^{i+j+1} b_i c_j
  \Pf \left( A^{\hat{i}, \hat{j}}_{\hat{i}, \hat{j}} \right)
  +
  \sum_{1 \le i < j \le n} (-1)^{i+j} b_i c_j
  \Pf \left( A^{\hat{i}, \hat{j}}_{\hat{i}, \hat{j}} \right)
  +
  d \Pf(A)
 \right).
\endalign
$$
By comparing this expression with the expansion of the second 
Pfaffian in (\EG) along the last two columns,
$$
\align
&\Pf \pmatrix A       & b & c \\
              -{}^t b & 0 & -d \\
              -{}^t c & d & 0 \endpmatrix
\\
&=
\sum_{1 \le i < j \le n} (-1)^{j-1} c_j \cdot (-1)^{i-1} b_i
 \Pf \left( A_{\hat{i}, \hat{j}}^{\hat{i}, \hat{j}} \right)
+
\sum_{1 \le j < i \le n} (-1)^{j-1} c_j \cdot (-1)^{i-2} b_i
 \Pf \left( A_{\hat{i}, \hat{j}}^{\hat{i}, \hat{j}} \right) \\
&\quad + (-d) \Pf(A).
\endalign
$$
we obtain the desired formula.

Next suppose that $n$ is odd.
Then $\det A = 0$.
By induction hypothesis, we see that, if $i < j$,
$$
\align
\det \left( A_{\hat{i}}^{\hat{j}} \right) 
&=
(-1)^{(n-i-1)+(n-j)} \det
 \pmatrix
       &      &     & a_{1i} \\
       & A_{\hat{i},\hat{j}}^{\hat{i},\hat{j}} && \vdots \\
       &      &     & a_{ni} \\
a_{j1} & \cdots & a_{jn} & a_{ji}
\endpmatrix
\\
&=
(-1)^{i+j-1}
\Pf
  \pmatrix
       &      &     & a_{1i} \\
       & A_{\hat{i},\hat{j}}^{\hat{i},\hat{j}} && \vdots \\
       &      &     & a_{ni} \\
a_{i1} & \cdots & a_{in} & 0
\endpmatrix
\Pf
  \pmatrix
       &      &     & a_{1j} \\
       & A_{\hat{i},\hat{j}}^{\hat{i},\hat{j}} && \vdots \\
       &      &     & a_{nj} \\
a_{j1} & \cdots & a_{jn} & 0
\endpmatrix
\\
&=
(-1)^{i+j-1}
 \cdot
 (-1)^{n-i-1} \Pf \left( A_{\hat{i}}^{\hat{i}} \right)
 \cdot
 (-1)^{n-j} \Pf \left( A_{\hat{j}}^{\hat{j}} \right)
\\
&=
 \Pf \left( A_{\hat{i}}^{\hat{i}} \right)
 \Pf \left( A_{\hat{j}}^{\hat{j}} \right).
\endalign
$$
Similarly, for $i > j$, we have
$$
\det \left( A_{\hat{i}}^{\hat{j}} \right)
=
 \Pf \left( A_{\hat{i}}^{\hat{i}} \right)
 \Pf \left( A_{\hat{j}}^{\hat{j}} \right).
$$
Comparison of this with the expansion of the Pfaffians in (\EH) along the last columns 
 completes the proof.\quad \quad \qed
\enddemo

Now we apply this lemma.
If $b$ is even, then we have
$$
\align
&\det
\pmatrix
N(X_n, X_n) & M_R(X_n) & M_{\{ n-b+(a+b)/2 \}}(X_n) \\
- \trans M_R(X_n) & 0 & 0 \\
N(x_{n+1}, X_n) & M_R(x_{n+1}) & x_{n+1}^{n-b+(a+b)/2}
\endpmatrix
\\
&\quad =
-
\Pf \pmatrix
N(X_n, X_n) & M_R(X_n) \\
- \trans M_R(X_n) & 0
\endpmatrix
\\
&\qquad\times
\Pf
\pmatrix
N(X_n, X_n) & M_R(X_n) & M_{\{ n-b+(a+b)/2 \}}(X_n) & N(X_n, x_{n+1}) \\
- \trans M_R(X_n) & 0 & 0 & -\trans M_R(x_{n+1}) \\
- \trans M_{\{ n-b+(a+b)/2 \}}(X_n) & 0 & 0 & -x_{n+1}^{n-b+(a+b)/2} \\
N(x_{n+1}, X_n) & M_R(x_{n+1}) & x_{n+1}^{n-b+(a+b)/2} & 0
\endpmatrix
\\
&\quad=
-
\Pf \pmatrix
N(X_n, X_n) & M_R(X_n) \\
- \trans M_R(X_n) & 0
\endpmatrix
\\
&\qquad\times
(-1)^{n-b+1}
\Pf \pmatrix
N(X_{n+1}, X_{n+1}) & M_Q(X_{n+1}) \\
- \trans M_Q(X_{n+1}) & 0
\endpmatrix
 \tag{\EC}
\endalign
$$
If $b$ is odd, then we have
$$
\align
&\det
\pmatrix
N(X_n, X_n) & M_R(X_n) & M_{\{ n-b+(a+b)/2 \}}(X_n) \\
- \trans M_R(X_n) & 0 & 0 \\
N(x_{n+1}, X_n) & M_R(x_{n+1}) & x_{n+1}^{n-b+(a+b)/2}
\endpmatrix
\\
&\quad =
\Pf \pmatrix
N(X_n, X_n) & M_R(X_n) & M_{\{ n-b+(a+b)/2 \}}(X_n) \\
- \trans M_R(X_n) & 0 & 0 \\
- \trans M_{\{ n-b+(a+b)/2 \}}(X_n) & 0
\endpmatrix
\\
&\qquad\times
\Pf
\pmatrix
N(X_n, X_n) & M_R(X_n) & N(X_n, x_{n+1}) \\
- \trans M_R(X_n) & 0 & 0 & -\trans M_R(x_{n+1}) \\
N(x_{n+1}, X_n) & M_R(x_{n+1}) & x_{n+1}^{n-b+(a+b)/2} & 0
\endpmatrix
\\
&\quad=
\Pf \pmatrix
N(X_n, X_n) & M_Q(X_n) \\
- \trans M_Q(X_n) & 0
\endpmatrix
\\
&\qquad\times
(-1)^{n-b}
\Pf \pmatrix
N(X_{n+1}, X_{n+1}) & M_R(X_{n+1}) \\
- \trans M_R(X_{n+1}) & 0
\endpmatrix
 \tag{\ED}
\endalign
$$

The four Pfaffians in (\EC) and (\ED) are evaluated by the following Lemma.

\proclaim{Lemma~\TJ}
{\rm(1)} If $a$ and $b$ is even, then
$$
\align
&\Pf \pmatrix
N(X_n, X_n) & M_R(X_n) \\
- \trans M_R(X_n) & 0
\endpmatrix
\\
&\quad=
 \frac{(-1)^{b(n-b)+(n-b)(n-b-1)/2}}{\Delta(X_n)} \\
&\qquad \times
  \det \left( M_{[0, n-b/2-1]}(X_n) \quad \vdots \quad M_{[n+a/2-b/2, n+a/2-1]}(X_n) \right) \\
&\qquad \times
  \det \left( M_{[0, n-b/2-1]}(X_n) \quad \vdots \quad M_{[n+a/2-b/2+1, n+a/2]}(X_n) \right),
\\
&\Pf \pmatrix
N(X_{n+1}, X_{n+1}) & M_Q(X_{n+1}) \\
- \trans M_Q(X_{n+1}) & 0
\endpmatrix
\\
&\quad=
 \frac{(-1)^{b(n-b)+b/2+(n-b+1)(n-b)/2}}{\Delta(X_{n+1})} \\
&\qquad \times
  \det \left( M_{[0, n-b/2]}(X_{n+1}) \quad \vdots \quad M_{[n+a/2-b/2+1, n+a/2]}(X_{n+1})
  \right) \\
&\qquad \times
  \det \left( M_{[0, n-b/2-1]}(X_{n+1}) \quad \vdots \quad M_{[n+a/2-b/2, n+a/2]}(X_{n+1})
 \right).
\endalign
$$
{\rm(2)} If $a$ and $b$ is odd, then
$$
\align
&\Pf \pmatrix
N(X_n, X_n) & M_Q(X_n) \\
- \trans M_Q(X_n) & 0
\endpmatrix
\\
&\quad=
 \frac{(-1)^{(b-1)(n-b)+(b-1)/2+(n-b+1)(n-b)/2}}{\Delta(X_n)} \\
&\qquad \times
  \det \left( M_{[0, n-(b-1)/2-1]}(X_n) \quad \vdots \quad
              M_{[n+(a+1)/2-(b-1)/2, n+(a+1)/2-1]}(X_n) \right) \\
&\qquad \times
  \det \left( M_{[0, n-(b+1)/2-1]}(X_n) \quad \vdots \quad
              M_{[n+(a+1)/2-(b+1)/2, n+(a+1)/2-1]}(X_n) \right),
\\
&\Pf \pmatrix
N(X_{n+1}, X_{n+1}) & M_R(X_{n+1}) \\
- \trans M_R(X_{n+1}) & 0
\endpmatrix
\\
&\quad=
 \frac{(-1)^{(b+1)(n-b)+(n-b)(n-b-1)/2}}{\Delta(X_{n+1})} \\
&\qquad \times
  \det \left( M_{[0, n-(b+1)/2]}(X_{n+1}) \quad \vdots \quad
              M_{[n+(a-1)/2-(b+1)/2+1, n+(a-1)/2]}(X_{n+1}) \right) \\
&\qquad \times
  \det \left( M_{[0, n-(b+1)/2]}(X_{n+1}) \quad \vdots \quad
              M_{[n+(a+1)/2-(b+1)/2+1, n+(a+1)/2]}(X_{n+1}) \right).
\endalign
$$
\endproclaim

\demo{Proof}
We use Theorem~\TG\ and the following decomposition of the product of 
two Schur functions of rectangular shape \cite{\OkadAI, Theorem~2.4}:
Let $m\le n$, then
$$
s_{(s^m)}(X_N) \cdot s_{(t^n)}(X_N)
 = \sum_{\lambda} s_{\lambda}(X_N),
$$
where the sum is taken over all partitions $\lambda$ with
 length $\le m+n$ such that
$$
\gather
\lambda_i + \lambda_{m+n-i+1} = s+t, \quad i=1, \dots, m, \\
\lambda_m \ge \max (s,t), \\
\lambda_{m+1}= \dots = \lambda_n = t.
\endgather
$$

Apply the minor summation formula in Theorem~\TG\ to
$$
G = M_P(X_n), \quad
H = M_R(X_n),
$$
and the skew-symmetric matrix $A = (a_{ij})_{i,j \in \Gamma}$ with nonzero
 entries
$$
a_{i,a+b-i} = \cases
 \hphantom{-}1  &\text{if $0 \le i \le (a+b)/2-1$} \\
 -1 &\text{if $(a+b)/2+1 \le i \le a+b$.}
\endcases
$$
Then, by the same argument as in the proof of \cite{\OkadAI, Theorem~2.4},
 we obtain the first formula in item (1) of the Theorem.
The other formulas are obtained by applying Theorem~\TG\ to the above
 skew-symmetric matrix $A$ and the matrices

$$
(G,H) = (M_P(X_{n+1}), M_Q(X_{n+1})),
\quad
(M_P(X_n), M_Q(X_n)),
\quad
(M_P(X_{n+1}), M_R(X_{n+1})). \qed
$$
\enddemo

Now we are in the position to complete the proof of Theorem~\TC. 

\demo{Proof of Theorem~\TC}
Suppose that $a$ and $b$ are even.
Combining \thetag{\EA}, \thetag{\EB} and \thetag{\EC}, and using
(\EF), we have
$$
\align
&\sum_{(\lambda, \mu) \in \Cal{R}(a,b)} s_{\lambda}(X_{n+1}) s_{\mu}(X_n) \\
&=
\frac{(-1)^{bn+n-b}}{\Delta(X_n) \Delta(X_{n+1})} \cdot
(-1)^{bn+b(b-1)/2+n-b} \cdot (-1)^{n-b}
\frac{(-1)^{(n-b)^2+b/2}}{\Delta(X_n) \Delta(X_{n+1})} \\
&\qquad\times
  \det \left( M_{[0, n-b/2-1]}(X_n) \quad \vdots \quad M_{[n+a/2-b/2, n+a/2-1]}(X_n) \right) \\
&\qquad \times
  \det \left( M_{[0, n-b/2-1]}(X_n) \quad \vdots \quad M_{[n+a/2-b/2+1, n+a/2]}(X_n) \right) \\
&\qquad \times
  \det \left( M_{[0, n-b/2]}(X_{n+1}) \quad \vdots \quad M_{[n+a/2-b/2+1, n+a/2]}(X_{n+1})
  \right) \\
&\qquad \times
  \det \left( M_{[0, n-b/2-1]}(X_{n+1}) \quad \vdots \quad M_{[n+a/2-b/2, n+a/2]}(X_{n+1})
 \right) \\
&=
(-1)^{b^2/2}
 s_{\left( \left( \frac{a}{2} \right)^{b/2} \right)}(X_n)
 s_{\left( \left( \frac{a}{2}+1 \right)^{b/2} \right)}(X_n)
 s_{\left( \left( \frac{a}{2} \right)^{b/2} \right)}(X_{n+1})
 s_{\left( \left( \frac{a}{2} \right)^{b/2+1} \right)}(X_{n+1})
\\
&=
 s_{\left( \left( \frac{a}{2} \right)^{b/2} \right)}(X_n)
 s_{\left( \left( \frac{a}{2}+1 \right)^{b/2} \right)}(X_n)
 s_{\left( \left( \frac{a}{2} \right)^{b/2} \right)}(X_{n+1})
 s_{\left( \left( \frac{a}{2} \right)^{b/2+1} \right)}(X_{n+1}).
\endalign
$$

If $a$ and $b$ are odd, then it follows from \thetag{\EA},
\thetag{\EB}, \thetag{\ED}, and (\EF) that

\vbox{
$$
\align
&\sum_{(\lambda, \mu) \in \Cal{R}(a,b)} s_{\lambda}(X_{n+1}) s_{\mu}(X_n) \\
&=
\frac{(-1)^{bn+n-b}}{\Delta(X_n) \Delta(X_{n+1})} \cdot
(-1)^{bn+b(b-1)/2+n-b} \cdot (-1)^{n-b}
\frac{(-1)^{2b(n-b)+(n-b)^2+(b-1)/2}}{\Delta(X_n) \Delta(X_{n+1})} \\
&\qquad \times
  \det \left( M_{[0, n-(b-1)/2-1]}(X_n) \quad \vdots \quad
              M_{[n+(a+1)/2-(b-1)/2, n+(a+1)/2-1]}(X_n) \right) \\
&\qquad \times
  \det \left( M_{[0, n-(b+1)/2-1]}(X_n) \quad \vdots \quad
              M_{[n+(a+1)/2-(b+1)/2, n+(a+1)/2-1]}(X_n) \right) \\
&\qquad \times
  \det \left( M_{[0, n-(b+1)/2]}(X_{n+1}) \quad \vdots \quad
              M_{[n+(a-1)/2-(b+1)/2+1, n+(a-1)/2]}(X_{n+1}) \right) \\
&\qquad \times
  \det \left( M_{[0, n-(b+1)/2]}(X_{n+1}) \quad \vdots \quad
              M_{[n+(a+1)/2-(b+1)/2+1, n+(a+1)/2]}(X_{n+1}) \right) \\
&=
(-1)^{(b^2-1)/2}
 s_{\left( \left( \frac{a+1}{2} \right)^{(b-1)/2} \right)}(X_n)
 s_{\left( \left( \frac{a+1}{2} \right)^{(b+1)/2} \right)}(X_n)\\
&\hskip4cm \times
 s_{\left( \left( \frac{a-1}{2} \right)^{(b+1)/2} \right)}(X_{n+1})
 s_{\left( \left( \frac{a+1}{2} \right)^{(b+1)/2} \right)}(X_{n+1}) \\
&=
 s_{\left( \left( \frac{a+1}{2} \right)^{(b-1)/2} \right)}(X_n)
 s_{\left( \left( \frac{a+1}{2} \right)^{(b+1)/2} \right)}(X_n)
 s_{\left( \left( \frac{a-1}{2} \right)^{(b+1)/2} \right)}(X_{n+1})
 s_{\left( \left( \frac{a+1}{2} \right)^{(b+1)/2} \right)}(X_{n+1}).
\endalign
$$
\line{\hfil\hbox{\qed}\quad \quad }
}
\enddemo


\Refs

\ref\no \CiucAB\by M.    Ciucu \yr 1997
\paper Enumeration of perfect matchings in graphs with reflective symmetry
\jour J. Combin\. Theory Ser.~A\vol 77 \pages 67-97\endref

\ref\no \CiucAH\by M.    Ciucu \yr \paper Enumeration of lozenge tilings 
of punctured hexagons\jour preprint\vol \pages \endref

\ref\no \FuKrAA\by M.    Fulmek and C. Krattenthaler \yr 1997 \paper 
Lattice path proofs for determinant formulas for symplectic and orthogonal 
characters\jour J. Combin\. Theory Ser.~A\vol 77\pages 3--50\endref

\ref\no \GeViAB\by I. M. Gessel and X. Viennot \yr 1989 
\paper Determinants, paths, and plane partitions 
\paperinfo preprint, 1989\endref

\ref\no \IsWaAA\by M.    Ishikawa and M. Wakayama \yr 1995 
\paper Minor summation formula for pfaffians
\jour Linear and Multilinear Algebra\vol 39
\pages 285--305\endref

\ref\no \IsWaAB\by M.    Ishikawa and M. Wakayama \yr 1995 
\paper Minor summation formula of pfaffians and Schur function identities
\jour Proc\. Japan Acad\. Ser.~A\vol 71
\pages 54--57\endref

\ref\no \IsOWAA\by M.    Ishikawa, S. Okada and M. Wakayama \yr 1996 
\paper Applications of minor summation formulas I, Littlewood's formulas
\jour J.~Algebra \vol 183 
\pages 193--216\endref

\ref\no \KupeAA\by G.    Kuperberg \yr 1994 
\paper Symmetries of plane partitions and the permanent determinant method 
\jour J.~Combin\. Theory Ser\. A \vol 68 
\pages 115--151\endref

\ref\no \MacdAC\by I. G. Macdonald \yr 1995 \book Symmetric Functions and 
Hall Polynomials \bookinfo second edition\publ Oxford University 
Press\publaddr New York/Lon\-don\endref

\ref\no \MacMAA\by P. A. MacMahon 
\book Combinatory Analysis 
\bookinfo vol.~2
\publ Cambridge University Press, 1916; reprinted by Chelsea, New York, 1960 \endref

\ref\no \OkadAI\by S.    Okada 
\paper Applications of minor summation formulas to rectangular-shaped representations of classical groups
\jour J. Algebra \toappear\vol 
\pages \endref

\ref\no \PropAA\by J.    Propp 
\paper Twenty open problems on enumeration of matchings
\jour preprint\vol 
\pages \endref

\ref\no \SagaAL\by B. E. Sagan \yr 1991 \book The symmetric group\publ 
Wadsworth \& Brooks/Cole\publaddr Pacific Grove, California\endref

\ref\no \StemAE\by J. R. Stembridge \yr 1990 
\paper Nonintersecting paths, pfaffians and plane partitions
\jour Adv\. in Math\.\vol 83
\pages 96---131\endref

\endRefs
\enddocument